\documentclass[12pt]{article}
\usepackage{bbm}
\usepackage{amsfonts}
\usepackage{pifont}
\usepackage{hyperref}
\usepackage{mathrsfs}
\usepackage{indentfirst}
\usepackage{amsmath}
\usepackage{amssymb}
\usepackage[amsmath, thmmarks]{ntheorem}
\usepackage{cite}
\usepackage{graphicx}
\usepackage{tikz}
\newtheorem{definition}{\bf Definition}[section]
\newtheorem{lemma}{\bf Lemma}[section]
\newtheorem{theorem}{\bf Theorem}[section]
\newtheorem{remark}{\bf Remark}[section]

\newtheorem{example}{\bf Example}[section]
\newtheorem{algorithm}{\bf Algorithm}[section]

\def\QEDopen{{\setlength{\fboxsep}{0pt}\setlength{\fboxrule}{0.2pt}\fbox{\rule[0pt]{0pt}{1.3ex}\rule[0pt]{1.3ex}{0pt}}}} 
\def\QED{\QEDopen}
\def\proof{{\bf Proof.} }
\def\endproof{\hspace*{\fill}~\QED\par\endtrivlist\unskip}

\begin{document}
\setcounter{page}{1}

\title{{\textbf{Finite atomic lattices and their monomial ideals }}\thanks {Supported by the National Natural Science
Foundation of China (No.11171242)}}
\author{Peng He\footnote{\emph{E-mail address}: 443966297@qq.com}, Xue-ping Wang\footnote{Corresponding author. xpwang1@hotmail.com; fax: +86-28-84761393}\\
\emph{College of Mathematics and Software
Science, Sichuan Normal University,}\\
\emph{Chengdu, Sichuan 610066, People's Republic of China}}
\newcommand{\pp}[2]{\frac{\partial #1}{\partial #2}}
\date{}
\maketitle

\begin{quote}
{\bf Abstract}

This paper primarily studies monomial ideals by their associated lcm-lattices. It first introduces notions of weak coordinatizations of finite atomic lattices which have weaker hypotheses than coordinatizations and shows the characterizations of all such weak coordinatizations. It then defines a finite super-atomic lattice in $\mathcal{L}(n)$, investigates the structures of $\mathcal{L}(n)$ by their super-atomic lattices and proposes an algorithm to calculate all the super-atomic lattices in $\mathcal{L}(n)$. It finally presents a specific labeling of finite atomic lattice and obtains the conditions that the specific labelings of  finite atomic lattices are the weak coordinatizations or the coordinatizations by using the terminology of super-atomic lattices. \\


\end{quote}

\section{Introduction}\label{intro}

Let $M$ be a monomial ideal in a polynomial ring $R=K[x_{1}, x_{2},\cdot\cdot\cdot, x_{n}]$ where $K$ is a field. We are interested in studying a minimal free resolution of $R/M$, and specifically understanding the maps in this resolution (see \cite{Dave, Clark09, Clark10, Eliahou, Velasco}). For a monomial ideal $M$, a minimal resolution is completely dependent on the information in the lcm-lattice of $M$, or LCM(M), which is the lattice of least common multiples of the minimal generators of $M$ partially ordered by divisibility. In 1999, Gasharov, Peeva, and Welker in \cite{Welker} expressed the multigraded Betti numbers of $R/M$ using the homology groups of certain open intervals in LCM(M). They further showed that the combinatorial type of minimal resolutions of a monomial ideal is determined by its $LCM$ lattice. In 2006, Phan in \cite{Phan} proved that all finite atomic lattices can be realized as the $LCM$ lattice of some monomial ideal $M$. He gave a construction which is motivated by the observation that for any coordinatization of an atomic lattice as a monomial ideal the set of lattice elements for which a given variable has a given degree bound is an order ideal. Essentially, he identified which order ideals are necessary and labels them with variables. In 2009, Mapes gave a generalization of the main construction in \cite{Phan} to describe all monomial ideals with a given $LCM$ lattice, i.e., she proved a statement as below (see \cite{Mapes09}, also \cite{Mapes13}).

Any labeling $\mathcal{M}$ of elements in a finite atomic lattice $P$ by monomials satisfying the following two conditions will yield a coordinatization of the lattice $P$.\\
$(A1)$ If $p\in \mbox{mi}(P)$ then $m_{p}\neq 1$ (i.e., all meet-irreducibles are labeled).\\
$(A2)$ If $\mbox{gcd}(m_{p}, m_{q})\neq 1$ for some $p,q \in P$ then $p$ and $q$ must be comparable (i.e., each variable only appears in monomials along one chain in $P$).

Mapes thought that it would be interesting to give an explicit formulation for when two coordinatizations are equivalent in this sense or to prove a version of the result above which has weaker hypotheses. This question has been inadvertently answered by Lukas Katth\"{a}n in \cite{Katthan} and separately by Maps and Piechnik in \cite{Mapes15} using different techniques. However, all of them do not give a general construction of the labeling $\mathcal{M}$ which does not satisfy the conditions (A1) and (A2) but $\mathcal{M}$ is a coordinatization.

On the other hand, the fact that the set of finite atomic lattices on $n$ ordered atoms, denoted by $\mathcal{L}({n})$, is itself a finite atomic lattice leads us to the question: what is the relationship between minimal resolutions of coordinatizations of lattices in $\mathcal{L}(n)$? The answer, due to a result in \cite{Welker}, is that the total Betti numbers are weakly monotonic along chains in $\mathcal{L}(n)$. This inspires us to understand the structure of $\mathcal{L}(n)$. In 2013, Mapes in \cite{Mapes13} proved that for any relation $P> Q$ in $\mathcal{L}(n)$ there exists a coordinatization of $Q$ producing a monomial ideal $M_{Q}$ and a deformation of exponents of $M_{Q}$ such that the lcm-lattice of the deformed ideal is $P$.

This paper will continue the topics on describing all monomial ideals by their $LCM$ lattices and understanding the structure of $\mathcal{L}(n)$, which is organized as follows. In Section 2, we give some preliminaries for convenience. In Section 3, we introduce notions of weak coordinatizations of finite atomic lattices and show their characterizations. In Section 4, we define a finite super-atomic lattice in $\mathcal{L}(n)$, investigate the structures of $\mathcal{L}(n)$ by their super-atomic lattices and propose an algorithm to calculate all the super-atomic lattices in $\mathcal{L}(n)$. In the end, we present a specific labeling of finite atomic lattice and obtain the conditions which are used to determine whether the specific labelings are the weak coordinatizations or the coordinatizations by terminology of super-atomic lattices.

\section{Preliminaries}
A poset is a structure $(P, \leq)$ where $P$ is a nonempty set and $\leq$ an ordering (reflexive, antisymmetric and transitive) relation on $P$.
We write $x\| y$ if $x\ngeq y$ and $y\ngeq x$, and we say that $x$ and $y$ are not comparable.
Contrarily, we write $x\nparallel y$ if $x\geq y$ or $y\geq x$, and we say $x$ and $y$ are comparable.
In addition, if $x< y$ and there is no element $z\in P$ such that $x< z< y$, then we say that $x$ is covered by $y$ (or $y$ covers $x$),
and we write $x\prec y$ (or $y\succ x$), see \cite{Crawley}.
\begin{definition}[\cite{Mapes13}]\label{sect2-2d}
\emph{A lattice is a poset $(P, \leq)$ satisfying the following properties:\\
(1) $P$ has a maximum element denoted by $1$.\\
(2) $P$ has a minimum element denoted by $0$.\\
(3) Every pair of elements $a$ and $b$ in $P$ has a join $a\vee b$ which is the least upper bound of the two elements.\\
(4) Every pair of elements $a$ and $b$ in $P$ has a meet $a\wedge b$ which is the greatest lower bound of the two elements.\\
If $P$ only satisfies conditions (2) and (4) then it is a meet-semilattice, and if $P$ only satisfies conditions (1) and (3) then it is a join-semilattice. Furthermore, if $P$ is a meet-semilattice with a unique maximal element then it is a lattice. Equivalently, if $P$ is a join-semilattice with a unique minimal element then it is a lattice.}
\end{definition}

We define an atom of a lattice $P$ to be an element $x\in P$ such that $x$ covers $0$. We denote the set of atoms in $P$ by atoms$(P)$ (see \cite{Crawley,Mapes13}). Let $A$ and $B$ be two sets. Then we denote that $A\setminus B=\{x\in A: x\notin B\}$, for convenience, if $B=\{b\}$ then we write $A\setminus B$ as $A\setminus b$.

\begin{definition}[\cite{Mapes13}]\label{sect2-3d}
\emph{If $P$ is a lattice and every element in $P\setminus 0$ is the join of atoms, then $P$ is an atomic lattice. Furthermore, if $P$ is finite, then it is a finite atomic lattice.}
\end{definition}

If $P$ is a lattice, then we define an element $x\in P$ to be meet-irreducible if $x\neq a\wedge b$ for any $a> x, b> x$.
We denote the set of meet-irreducible elements in $P$ by mi$(P)$. Given an element $x\in P$, an order ideal of $x$ is defined to be the set $\lfloor x\rfloor=\{a\in P: a\leq x\}$. Similarly, we define an order filter of $x$ to be $\lceil x\rceil=\{a\in P: x\leq a\}$ (see \cite{Crawley,Mapes13}).

\begin{lemma}[ Lemma 2.3 of \cite{Mapes13}]\label{sect2-1l}
Let $P$ be a finite atomic lattice. Every element $p\in P$ is the meet of all the meet-irreducible elements $l$ such that $l\geq p$.
\end{lemma}

It will be convenient to consider finite atomic lattices as sets of sets in the following way. Let $\mathcal{S}$ be a set of subsets of $\{1,\cdot \cdot \cdot,n\}$ with no duplicates, closed under intersections, and containing the entire set, the empty set, and the sets $\{i\}$ for all $1\leq i\leq n$. Then it is easy to see that $\mathcal{S}$ is a finite atomic lattice by ordering the sets in $\mathcal{S}$ by inclusion. Conversely, it is clear that any finite atomic lattice $P$ can be expressed in this way, simply by letting
$$\mathcal{S}_{P}=\{\sigma: \sigma=\mbox{supp}(p), p\in P\},$$
where $\mbox{supp}(p)=\{a_{i}: a_{i}\leq p, a_{i}\in \mbox{atoms}(P)\}$ (see \cite{Birkhoff,Buchi,Mapes13}).

\begin{definition}[\cite{Welker}]\label{sect2-4d}
\emph{The $LCM$ lattice, $LCM(M)$, of a monomial ideal $M$ is the set of least common multiples of minimal generators of $M$, partially ordered by divisibility.}
\end{definition}

\begin{example}\label{sect2-1e}
\emph{For the monomial ideal $M = (a^{2}cd, abd, abc)\subseteq k[a,b,c,d]$, the Hasse diagram of the $LCM$ lattice of $M$ is shown as Fig.1 (note the minimal element of the lattice has been left off, as will often be the case).}
\par\noindent\vskip50pt
 \begin{minipage}{11pc}

\setlength{\unitlength}{0.75pt}\begin{picture}(300,120)

\put(250,40){\circle{4}}\put(240,28){\makebox(0,0)[l]{\footnotesize $abd$}}

\put(350,40){\circle{4}}\put(360,28){\makebox(0,0)[r]{\footnotesize$ abc$}}

\put(150,40){\circle{4}}\put(140,28){\makebox(0,0)[l]{\footnotesize$a^{2}cd$}}

\put(300,90){\circle{4}}\put(330,100){\makebox(0,0)[r]{\footnotesize$abcd$}}

\put(250,140){\circle{4}}\put(260,150){\makebox(0,0)[r]{\footnotesize$a^{2}bcd$}}

  \put(251,41){\line(1,1){47}}
  \put(151,41){\line(1,1){98}}
  \put(349,41){\line(-1,1){47}}
  \put(299,91){\line(-1,1){47}}

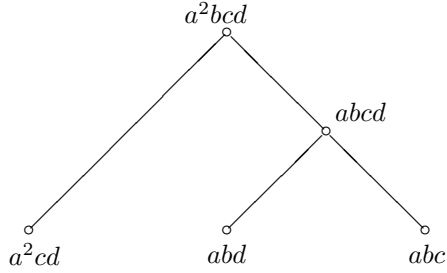
 \put(155,0){\emph{Fig.1. The lattice $LCM(M)$}}
 \end{picture}
 \end{minipage}
  \end{example}

One result in \cite{Welker} is that for monomial ideals all minimal resolutions are completely dependent on the information in the $LCM$ lattice. Specifically, one can compute multigraded Betti numbers using the $LCM$ lattice $LCM(M)$ and all ideals with a given $LCM$ lattice have isomorphic minimal free resolutions.

\begin{definition}[\cite{Mapes09}]\label{sect2-5d}
\emph{Define a labeling of a finite atomic lattice $P$ to be any assignment of non-trivial monomials $\mathcal{M}=\{m_{p_{1}},\cdot \cdot\cdot, m_{p_{t}}\}$ to some set of elements $p_{i}\in P$. It will be convenient to think of unlabeled elements as having the label 1. Define a monomial ideal $M_{P,\mathcal{M}}$ to be the ideal generated by monomials
\begin{equation}\label{eqhe0}x(a)=\prod_{p\in\lceil a\rceil^{c}}m_{p}\end{equation}
for each $a\in \mbox{atoms}(P)$ where $\lceil a\rceil^{c}$ means taking the complement of $\lceil a\rceil$ in $P$. We say that the labeling $\mathcal{M}$ is a coordinatization if the lcm-lattice of $M_{P,\mathcal{M}}$ is isomorphic to $P$.}
\end{definition}

\begin{lemma}[Proposition 3.2.1 of \cite{Mapes09} and Theorem 3.2 of \cite{Mapes13},]\label{sect2-2l}
Any labeling $\mathcal{M}$ of elements in a finite atomic lattice $P$ by monomials satisfying the following two conditions will yield a coordinatization of the lattice $P$.\\
$(A1)$ If $p\in \emph{mi}(P)$ then $m_{p}\neq 1$ (i.e., all meet-irreducibles are labeled).\\
$(A2)$ If $\emph{gcd}(m_{p}, m_{q})\neq 1$ for some $p,q \in P$ then $p\nparallel q$ (i.e., each variable only appears in monomials along one chain in $P$).
\end{lemma}

Let $\mathcal{M}$ be a labeling with conditions (A1) and (A2) hold, and let $f: P\rightarrow LCM(M_{P,\mathcal{M}})$ be denoted by
\begin{equation}\label{eqhe}f(p)=\prod_{q\in \lceil p\rceil^{c}}m_{q}\end{equation} for each $p\in P$.
Then $f$ is an isomorphism from $P$ to $LCM(M_{P,\mathcal{M}})$.

\begin{lemma}[Lemma 3.3 of \cite{Mapes13}]\label{sect2-3l}
If $p\in \lceil q\rceil^{c}$ for some $p, q \in P$ where $P$ is a finite atomic lattice, then $\lfloor p\rfloor \subseteq \lceil q\rceil^{c}$.
\end{lemma}

Let $M$ be a monomial ideal with $n$ generators and let $P_{M}$ be its lcm-lattice. For notational purposes, denote $P_{M}$ as the set consisting of elements denoted $\bar{p}$ which represent the monomials occurring in $P_{M}$. Now, define an abstract finite atomic lattice $P$ where the elements in $P$ are formal symbols $p$ satisfying the relations $p< p^{'}$ if and only if $\bar{p}< \bar{p}^{'}$ in $P_{M}$. In other words, $P$ is the finite atomic lattice isomorphic to $P_{M}$ obtained by simply forgetting the data of the monomials in $P_{M}$. Define a labeling of $P$ in the following way, let $\mathcal{D}$ be the set consisting of monomials $m_{p}$ for each $p\in P$ defined by
\begin{equation}\label{eq1} m_{p}=\frac{\mbox{gcd}\{\bar{t}: t> p\}}{\bar{p}},\end{equation} where by convention $\mbox{gcd}\{\bar{t}: t> p\}$ for $p=1$ is defined to be $\bar{1}$. Note that $m_{p}$ is a monomial since clearly $\bar{p}$ divides $\bar{t}$ for all $t> p$.

\begin{lemma}[Proposition 3.6 of \cite{Mapes13}]\label{sect2-4l}
Given $M$ a monomial ideal with lcm-lattice $P_{M}$. If $P$ is an abstract finite atomic lattice where $P$ is isomorphic to $P_{M}$ as lattices then the labeling $\mathcal{D}$ of $P$ as defined by (\ref{eq1}) is a coordinatization and the resulting monomial ideal $M_{P,\mathcal{D}}=M$.
\end{lemma}

Although Lemma \ref{sect2-4l} shows that the labeling $\mathcal{D}$ of $P$ as defined by (\ref{eq1}) is a coordinatization, the following theorem will further verify that the labeling $\mathcal{D}$ induced by (\ref{eq1}) is the same as $\mathcal{M}$ if $\mathcal{M}$ satisfies the conditions of Lemma \ref{sect2-2l}.

Customarily, we denote that $\mbox{lcm }\emptyset=1$ and $\mbox{gcd }\emptyset=1$.
\begin{theorem}\label{sect3-2t}
Let $\mathcal{M}=\{m_p: p\in P\}$ be a labeling of a finite atomic lattice $P$ satisfying the conditions of Lemma \ref{sect2-2l}, and let $M=M_{P,\mathcal{M}}$
and for each $p\in P$, $\bar{p}=f(p)$ where $f(p)$ is defined by (\ref{eqhe}). Then the labeling $\mathcal{D}=\{m^{'}_p: p\in P\}$ of $P$ as defined by (\ref{eq1})
satisfies $m^{'}_{p}=m_{p}$ for each $p\in P$.
\end{theorem}
\proof
Suppose that $P$ has $n$ atoms. We first note that $\bar{p}=f(p)=\prod_{q\in \lceil p\rceil^{c}}m_{q}$ for all $p\in P$. Thus the formula (\ref{eq1}) implies that
\begin{eqnarray*}m_{p}^{'}&=&\frac{\mbox{gcd}\{\prod_{q\in \lceil t\rceil^{c}} m_{q}: t> p\}}{\prod_{q\in \lceil p\rceil^{c}}m_{q}}\\
&=&\frac{\prod_{q\in \lceil p\rceil^{c}}m_{q}\ast \mbox{gcd}\{\prod_{q\in \lceil t\rceil^{c}\setminus\lceil p\rceil^{c}} m_{q}: t> p\}}{\prod_{q\in \lceil p\rceil^{c}}m_{q}}\\
&=&\mbox{gcd}\{\prod_{q\in \lceil t\rceil^{c}\setminus\lceil p\rceil^{c}}m_{q}: t> p\}.\end{eqnarray*}
Note that if $a\geq b$, then $\lceil a\rceil^{c} \supseteq \lceil b\rceil^{c}$, which means
$$\prod_{q\in \lceil b\rceil^{c}\setminus\lceil p\rceil^{c}}m_{q}\mid \prod_{q\in \lceil a\rceil^{c}\setminus\lceil p\rceil^{c}}m_{q}.$$
Thus $$m_{p}^{'}=\mbox{gcd}\{\prod_{q\in \lceil t\rceil^{c}\setminus\lceil p\rceil^{c}}m_{q}: t> p\}
=\mbox{gcd}\{\prod_{q\in \lceil t\rceil^{c}\setminus\lceil p\rceil^{c}}m_{q}: t\succ p\}.$$
This follows that $m_{p}^{'}= m_p\ast \mbox{gcd}\{\prod_{q\in (\lceil t\rceil^{c}\setminus\lceil p\rceil^{c})\setminus p}m_{q}: t\succ p\}$
since $p\in \lceil t\rceil^{c}\setminus\lceil p\rceil^{c}$ for any $t\succ p$.

Therefore, in order to prove $m_{p}=m_{p}^{'}$ for all $p\in P$, we just need to show
$$\mbox{gcd}\{\prod_{q\in (\lceil t\rceil^{c}\setminus\lceil p\rceil^{c})\setminus p}m_{q}: t\succ p\}=1$$ as follows.

(a) If there is only one element $t\in P$ satisfying $t\succ p$,
then $(\lceil t\rceil^{c}\setminus\lceil p\rceil^{c})\setminus p=\emptyset$.
Otherwise, there exists an element $d\in P$ such that $d> p$ and $d\ngeqslant t$, where $d\ngeqslant t$ implies that $d< t$ or $d\| t$.
If $d< t$ then $t>d>p$, contrary to $t\succ p$. If $d\| t$ then we have an element $c\in P$ such that $d\geq c\succ p$ since $d> p$.
Thus $c = t$, and then $d\geq t$, a contradiction.
Therefore, $\mbox{gcd}\{\prod_{q\in (\lceil t\rceil^{c}\setminus\lceil p\rceil^{c})\setminus p}m_{q}: t\succ p\}=\mbox{gcd }\emptyset=1$.

(b) Suppose that there are $k$ elements $t_{1}, t_{2}, \cdots, t_k$ in $P$ such that $t_{i}\succ p$ for any $1 \leq i \leq k$ where $k\geq 2$.
If $\mbox{gcd}\{\prod_{q\in (\lceil t\rceil^{c}\setminus\lceil p\rceil^{c})\setminus p}m_{q}: t\succ p\}\neq 1$,
then there exists a variable $x_p$ such that $x_p\mid \mbox{gcd}\{\prod_{q\in (\lceil t\rceil^{c}\setminus\lceil p\rceil^{c})\setminus p}m_{q}: t\succ p\}$.
Therefore, we have an element $q_i> p$
and $q_i\ngeq t_i$ such that $x_p \mid m_{q_{i}}$ for each $1\leq i\leq k$. By (A2), $\{q_1, q_2, \cdots, q_k\}$ lies in a chain in $P$.
Hence, there exists an element $1\leq r\leq k$ such that $\{q_1, q_2, \cdots, q_k, t_r\}$ be a chain,
and then for all $1\leq j\leq k$ we have $q_j\geq t_r$ since $q_j> p$ and $t_r\succ p$. Thus $q_r \geq t_r$, a contradiction.
Therefore, $\mbox{gcd}\{\prod_{q\in (\lceil t\rceil^{c}\setminus\lceil p\rceil^{c})\setminus p}m_{q}: t\succ p\}= 1$.
\endproof

\section{Weak coordinatizations}

One of the main results in \cite{Phan} is that every finite atomic lattice is in fact the lcm-lattice of a monomial ideal. In 2009, Mapes in \cite{Mapes09}
introduced a definition of coordinatization. Moreover, she proved that there are some specific constructions which produce a monomial ideal whose
lcm-lattice has a given lattice structure, i.e., Lemma \ref{sect2-2l} (see also \cite{Mapes13}).
Mapes thought that it would be interesting to give an explicit formulation for when two coordinatizations are equivalent in this sense or to
prove a version of Lemma \ref{sect2-2l} which has weaker hypotheses.

In this section, we shall introduce the notion of a weak coordinatization which has weaker hypotheses than Definition \ref{sect2-5d}, and show
a sufficient condition which yields a weak coordinatization.












Let $P$ be a finite atomic lattice and $p\in P$. Define $$B_p=\{T\subseteq \mbox{supp}(p): \bigvee_{b\in T}b=p\}.$$

\begin{definition}\label{sect3-1d}
\emph{
Let $\mathcal{M}$ be a labeling of a finite atomic lattice $P$. Define a monomial ideal $I_{P,\mathcal{M}}$ to be the ideal generated by monomials
\begin{equation}\label{equ1.2}\triangle(a)=\mbox{gcd}\{\mbox{lcm} \{x(b): b\in T\}: T\in \bigcup_{p\geq a}B_p\} \end{equation}
for each $a\in \mbox{atoms}(P)$. We say that the labeling $\mathcal{M}$ is a weak coordinatization if the lcm-lattice of $I_{P,\mathcal{M}}$ is isomorphic to $P$.}
\end{definition}

We first have the following lemma.
\begin{lemma}\label{sect3-he}
A labeling $\mathcal{M}$ is a coordinatization of a finite atomic lattice $P$ if and only if it is a weak coordinatization
and $\triangle(a)= x(a)$ for all $a\in \emph{atoms}(P)$.
\end{lemma}
\proof
By Definition\ref{sect3-1d}, the sufficiency is clear. Now, we prove the necessity.

First, for all $a\in \mbox{atoms}(P)$, as $\{a\}\in \bigcup \limits_{p\geq a}B_p$, equation (\ref{equ1.2}) implies $\triangle(a)\mid x(a)$.

Secondly, since $\mathcal{M}$ is a coordinatization, the map
$$g: P\rightarrow LCM(M_{P,\mathcal{M}}) \mbox{ with }g(a)=x(a)$$
for all $a\in \mbox{atoms}(P)$ is an isomorphism.
Thus, for any $p\in P$ and any $T\in B_p$, $$g(p)= \mbox{lcm}\{x(b): b\in \mbox{supp}(p)\} = \mbox{lcm}\{x(b): b\in T\}.$$

Finally, suppose that $a\in \mbox{atoms}(P)$. Let $p\in P$ and $a\leq p$. Clearly, $a\in \mbox{supp}(p)$ and then
$$g(a)=x(a)\mid g(p)=\mbox{lcm}\{x(b): b\in T\} \mbox{ for any }T\in B_p.$$
So that $x(a)\mid\mbox{lcm}\{x(b): b\in T\}$ for any $T\in \bigcup \limits_{p\geq a}B_p$.
Further, by (\ref{equ1.2}), $x(a)\mid \triangle(a)$.

Therefore, $\triangle(a)=x(a)$ for all $a\in \mbox{atoms}(P)$, which together with the fact that $\mathcal{M}$ is
a coordinatization of $P$ yields that $\mathcal{M}$ is a weak coordinatization of $P$.
\endproof

Notice that a weak coordinatization of a finite atomic lattice $P$ needs not to be a coordinatizaton. For instance,
let $P$ be the finite atomic lattice with a labeling as Fig.2. Then by Definitions \ref{sect2-5d} and \ref{sect3-1d},
 $$M_{P,\mathcal{M}}=(b^{2}c^{2}d^{2}e^{2}, acd^{2}e^{2}, a^{2}b^{2}d^{2}e^{2}, a^{2}b^{3}c^{2}e, a^{2}b^{3}c^{2}d),$$ $$I_{P,\mathcal{M}}=(b^{2}c^{2}d^{2}e^{2}, acd^{2}e^{2}, a^{2}b^{2}d^{2}e^{2}, a^{2}b^{2}c^{2}e, a^{2}b^{2}c^{2}d).$$
Then it is obvious that the lattice $LCM(I_{P,\mathcal{M}})$ shown as Fig.3 is isomorphic to $P$. Further, the labeling $\mathcal{M}$ is a weak coordinatizaton of $P$.
On the other hand, the lattice $LCM(M_{P,\mathcal{M}})$ shown as Fig.4 is not isomorphic to $P$, it follows that $\mathcal{M}$ is not a coordinatizaton of $P$.
 \par\noindent\vskip50pt
\begin{minipage}{11pc}
  \setlength{\unitlength}{0.75pt}\begin{picture}(160,30)
\put(-10,20){\circle{4}}\put(-20,15){\makebox(0,0)[l]{\footnotesize $a$}}
\put(30,20){\circle{4}}\put(20,15){\makebox(0,0)[l]{\footnotesize $b$}}
\put(70,20){\circle{4}}\put(80,15){\makebox(0,0)[r]{\footnotesize $c$}}
\put(110,20){\circle{4}}\put(120,15){\makebox(0,0)[r]{\footnotesize $d$}}
\put(150,20){\circle{4}}\put(158,15){\makebox(0,0)[r]{\footnotesize $e$}}
\put(130,40){\circle{4}}\put(147,45){\makebox(0,0)[r]{\footnotesize $de$}}
\put(10,40){\circle{4}}\put(-5,45){\makebox(0,0)[l]{\footnotesize $ab$}}
\put(50,40){\circle{4}}\put(65,45){\makebox(0,0)[r]{\footnotesize $bc$}}
\put(30,60){\circle{4}}
 \put(-9,21){\line(1,1){18}}
  \put(29,21){\line(-1,1){18}}
  \put(31,21){\line(1,1){18}}
  \put(69,21){\line(-1,1){18}}
  \put(111,21){\line(1,1){18}}
  \put(149,21){\line(-1,1){18}}
  \put(11,41){\line(1,1){18}}
  \put(49,41){\line(-1,1){18}}
  \put(129,41){\line(-5,1){97}}
 \put(-13,-10){Fig.2. $P$ with a labeling}
\put(190,20){\circle{4}}
\put(230,20){\circle{4}}
\put(270,20){\circle{4}}
\put(310,20){\circle{4}}
\put(350,20){\circle{4}}
\put(330,40){\circle{4}}
\put(210,40){\circle{4}}
\put(250,40){\circle{4}}
\put(230,60){\circle{4}}
  \put(191,21){\line(1,1){18}}
  \put(229,21){\line(-1,1){18}}
  \put(231,21){\line(1,1){18}}
  \put(269,21){\line(-1,1){18}}
  \put(311,21){\line(1,1){18}}
  \put(349,21){\line(-1,1){18}}
  \put(211,41){\line(1,1){18}}
  \put(249,41){\line(-1,1){18}}
  \put(329,41){\line(-5,1){97}}
 \put(201,-10){Fig.3. $LCM(I_{P,\mathcal{M}})$}
\put(390,20){\circle{4}}
\put(430,20){\circle{4}}
\put(470,20){\circle{4}}
\put(510,20){\circle{4}}
\put(550,20){\circle{4}}
\put(530,40){\circle{4}}
\put(490,80){\circle{4}}
\put(410,40){\circle{4}}
\put(450,40){\circle{4}}
\put(430,60){\circle{4}}
  \put(391,21){\line(1,1){18}}
  \put(429,21){\line(-1,1){18}}
  \put(431,21){\line(1,1){18}}
  \put(469,21){\line(-1,1){18}}
  \put(511,21){\line(1,1){18}}
  \put(549,21){\line(-1,1){18}}
  \put(411,41){\line(1,1){18}}
  \put(449,41){\line(-1,1){18}}
  \put(529,41){\line(-1,1){37.5}}
  \put(431,61){\line(3,1){58}}
 \put(406,-10){Fig.4. $LCM(M_{P,\mathcal{M}})$}
 \end{picture}
 \end{minipage}\\

\begin{lemma}\label{sect3-1l}
Let $\mathcal{M}$ be a labeling of a finite atomic lattice $P$ and $p\in P$. For each $R\in B_p$, if $b\in \emph{supp}(p)\setminus R$
then $\triangle(b)\mid \emph{lcm}\{\triangle(r): r\in R\}$.
\end{lemma}
\proof
Suppose $\triangle(b)\nmid \mbox{lcm}\{\triangle(r): r\in R\}$. Then there is a monomial $x^{u_{b}}$ such that $x^{u_{b}}\nmid \mbox{lcm}\{\triangle(r): r\in R\}$
where $x^{u_{b}}$ is the highest power of $x$ dividing $\triangle(b)$.
Let $$S=\{a\in R: x^{u_{b}}\mid x(a)\}$$ and $x^{u_{a}}$ be the highest power of $x$ dividing $\triangle(a)$ for each $a\in S$.
Then $u_{a}< u_{b}$ since $x^{u_{b}}\nmid \mbox{lcm}\{\triangle(r): r\in R\}$. Moreover, it follows from formula (\ref{equ1.2}) that for any $a\in S$ there exists an element $q_a\in P$ with $q_a\geq  a$ and a set $T_a\in B_{q_a}$ such that $x^{u_{a}}$ is the highest power of $x$ dividing $\mbox{lcm}\{x(t): t\in T_a\}$. Thus \begin{equation}\label{00000}x^{u_b}\nmid \mbox{lcm}\{x(t): t\in T_a\}\end{equation} for each $a\in S$ since $u_a < u_b$.

Next, let $C=\bigcup_{a\in S}T_a\bigcup (R\setminus S)$. Clearly, we have
$$\bigvee_{c\in C} c=\bigvee_{a\in S} (\bigvee T_a)\vee\bigvee(R\setminus S)=\bigvee_{a\in S}q_a\vee\bigvee(R\setminus S)\geq  \bigvee_{a\in S}a \vee\bigvee(R\setminus S)=p \geq b$$ and $C\in B_{\bigvee_{c\in C}c}$. Using (\ref{equ1.2}), we have $\triangle(b)\mid \mbox{lcm}\{x(c): c\in C\}$.
Thus \begin{equation}\label{00001}x^{u_{b}}\mid \mbox{lcm}\{x(c): c\in C\}.\end{equation} However, from (\ref{00000})
we know that if $c\in \bigcup_{a\in S}T_a$ then $x^{u_{b}}\nmid x(c)$. Moreover, if $c\in R\setminus S$ then $x^{u_{b}}\nmid x(c)$
by the construction of $S$. Hence, $x^{u_{b}}\nmid \mbox{lcm}\{x(c): c\in C\}$, contrary to (\ref{00001}).
Therefore, $\triangle(b)\mid \mbox{lcm}\{\triangle(r): r\in R\}$.
\endproof

\begin{lemma}\label{sect3-0l}
Let $\mathcal{M}$ be a labeling of a finite atomic lattice $P$. For all $p, q\in P$, if $x_0\mid m_p$ and $x_0\mid m_q$
imply $p\nparallel q$ then $x_0 \nmid \frac{x(a)}{\emph{gcd}(\triangle(a), x(a))}$ for any $a\in \emph{atoms}(P)$.
\end{lemma}
\proof
Let $S=\{s\in P: x_0 \nmid m_s\}$ and $R= P\setminus S$. Suppose that $\overline{m_s}=x_s$ with $s\in S$ and $\overline{m_r}= x_{0}^{r}$ where $x_{0}^{r}$ is the
highest power of $x_0$ dividing $m_r$ with $r\in R$. Then from the hypotheses of Lemma \ref{sect3-0l}, the labeling $\overline{\mathcal{M}}=\{\overline{m_p}: p\in P\}$
satisfies the conditions of Lemma \ref{sect2-2l}. Thus $\overline{\mathcal{M}}$ is a coordinatization of $P$. Hence, by Lemma \ref{sect3-he}, $\overline{\mathcal{M}}$ is a
weak coordinatization of $P$ and \begin{equation}\label{11} \overline{x(a)}=\overline{\triangle(a)}\end{equation} for any atom $a \in \mbox{atoms}(P)$
where $\overline{x(a)} \in \overline{M}_{P,\mathcal{\overline{M}}}$ and $\overline{\triangle(a)}\in I_{P,\mathcal{\overline{M}}}$.

Now, assume that $x_{0}^{a_1}$ and $x_{0}^{\overline{a_1}}$ are highest powers of $x_0$ dividing $x(a)$ and $\overline{x(a)}$, respectively,
and $x_{0}^{a_2}$ and $x_{0}^{\overline{a_2}}$ are highest powers of $x_0$ dividing $\triangle(a)$ and $\overline{\triangle(a)}$, respectively.
By Definition \ref{sect2-5d}, we have $a_1= \overline{a_1}$, which together with equation (\ref{equ1.2}) implies that $a_2=\overline{a_2}$. Using (\ref{11}), we have $\overline{a_1}=\overline{a_2}$.  Therefore, $a_1 =a_2$, which means $x_0 \nmid \frac{x(a)}{\mbox{gcd}(\triangle(a), x(a))}$.
\endproof

\begin{theorem}\label{sect3-1t}
Any labeling $\mathcal{M}$ of elements in a finite atomic lattice $P$ by monomials satisfying the following two conditions will yield a weak coordinatization of the lattice $P$.\\
\emph{(C1)} If $p\in \emph{mi}(P)$ then $m_{p}\neq 1$.\\
\emph{(C2)} If $\emph{gcd}(m_{p}, m_{q})\neq 1$ for some $p, q \in P$ then either $p\nparallel q$, or
$$r_q(p)=\frac{m_{p}}{\emph{gcd}(m_{p}, m_{q})}\neq 1, r_p(q)=\frac{m_{q}}{\emph{gcd}(m_{p}, m_{q})}\neq 1$$
and if $x, y\in\{s\in P: \emph{gcd}(r_q(p), m_s)\neq 1\}$ or $x, y\in\{s\in P: \emph{gcd}(r_p(q), m_s)\neq 1\}$ then $x\nparallel y$.
\end{theorem}
\proof
The proof of Theorem \ref{sect3-1t} is made in several steps. Let $P^{'}$ be the lcm-lattice of $I_{P,\mathcal{M}}$. For $b\in P$, define $g: P\longrightarrow P^{'}$ to be the map such that
\begin{equation}\label{11113}g(b)= \mbox{lcm}\{\triangle(a_{i}): a_{i}\in \mbox{supp}(b)\}.\end{equation} Next, we shall show that $g$ is an isomorphism from $P$ to $P^{'}$. Note that $g$ is well-defined.\\
$ $\\
A. $\triangle (a)\nmid \triangle(b)$ and $\triangle(b)\nmid\triangle(a)$ for any $a, b\in \mbox{atoms}(P)$ with $a\neq b$.

By Lemma \ref{sect2-1l}, the condition $a\neq b$ yields that $\mbox{mi}(P)\cap \lceil a\rceil \neq\mbox{mi}(P)\cap \lceil b\rceil$.
Moreover, $a\| b$ since $a\neq b$ and $a, b\in \mbox{atoms}(P)$. Thus by Lemma \ref{sect2-1l} $$\mbox{mi}(P)\cap \lceil a\rceil \nsubseteq \mbox{mi}(P)\cap \lceil b\rceil
\mbox{ and }\mbox{mi}(P)\cap \lceil b\rceil \nsubseteq \mbox{mi}(P)\cap \lceil a\rceil.$$ Hence
$$\mbox{mi}(P)\cap \lceil a\rceil^{c} \nsubseteq \mbox{mi}(P)\cap \lceil b\rceil^{c}\mbox{ and }\mbox{mi}(P)\cap \lceil b\rceil^{c} \nsubseteq \mbox{mi}(P)\cap \lceil a\rceil^{c}.$$
Therefore, there exists at least one element \begin{equation}\label{qqq}q\in \mbox{mi}(P)\cap \lceil a\rceil^{c}\mbox{ but }q\notin \mbox{mi}(P)\cap \lceil b\rceil^{c}.\end{equation}

We shall prove the following statement.
\begin{equation}\label{11110}\mbox{There exists a variable }x_q \mid m_q \mbox{ such that for all }r\in P, x_q\mid m_r  \mbox{ implies that }q\nparallel r.\end{equation}

Indeed, since $q$ is meet-irreducible, condition (C1) yields that $m_{q}\neq 1$.
Let $y_{q}$ be a variable satisfying $y_{q} \mid m_{q}$. Then there are two cases.

Case (1). If for all $r\in P$, $y_q\mid m_r$ implies $q\nparallel r$, then clearly (\ref{11110}) is true.

Case (2). If there is $t\in P$ such that $y_q\mid m_t$ but $q\|t$, then $\mbox{gcd}(m_{t}, m_{q})\neq 1$. Thus $r_t(q)\neq 1$ by condition (C2). Let $x_q\mid r_t(q)$ and $C_{q}=\{u\in P: x_q\mid m_u\}$.
Then $q\in C_{q}$ and $C_{q}\subseteq\{s\in P: \mbox{gcd}(r_t(q), m_s)\neq 1\}$. Again, by condition ($C_2$), $x\nparallel y$ for any $x, y \in \{s\in P: \mbox{gcd}(r_t(q), m_s)\neq 1\}$, i.e., $\{s\in P: \mbox{gcd}(r_t(q), m_s)\neq 1\}$ is a chain in $P$. Thus the condition $C_{q}\subseteq\{s\in P: \mbox{gcd}(r_t(q), m_s)\neq 1\}$ means that $C_{q}$ is a chain in $P$.
Note that $x_q\mid m_q$. Therefore, by the construction of $C_{q}$, we have that for all $r\in P$,
$x_q \mid m_r$ implies that $q\nparallel r$, i.e., (\ref{11110}) is true.

In the view of Cases (1) and (2), (\ref{11110}) holds.

Now let $x_{q}$ be a variable of $m_{q}$ such that (\ref{11110}) holds and let $D_{q}=\{v\in P: x_q\mid m_v\}$. Then $q\in D_q$. Suppose that $p\in \lceil b\rceil^{c}$ satisfies $x_{q}\mid m_{p}$. Then $p\nparallel q$ by (\ref{11110}).
Note that $p\neq q$. Thus, either $q< p$ or $p< q$. If $q< p$ then $q\in \lfloor p\rfloor \subseteq \lceil b\rceil^{c}$ by $p\in \lceil b\rceil^{c}$
and Lemma \ref{sect2-3l}, contrary to (\ref{qqq}). So that $p< q$. Therefore, for all $p\in \lceil b\rceil^{c}$, if $x_{q}\mid m_{p}$ then $p< q$.
Further, from the construction of $D_q$, we know that if $z\in D_{q}\cap \lceil b\rceil^{c}$ then $z< q$.
Note that $q\in \lceil a\rceil^{c}$ by (\ref{qqq}). Thus $z< q\in \lceil a\rceil^{c}$, it follows from Lemma \ref{sect2-3l} that $z\in \lceil a\rceil^{c}$.
So, $D_{q}\cap \lceil b\rceil^{c}\subseteq D_{q}\cap\lceil a\rceil^{c}$. Note that $q\in \lceil a\rceil^{c}, q\in D_q$
and $q\notin D_{q}\cap \lceil b\rceil^{c}$. Therefore, \begin{equation}\label{chen1}D_{q}\cap \lceil b\rceil^{c}\subsetneq D_{q}\cap\lceil a\rceil^{c}.\end{equation}

Finally, let $x_{q}^{s_a}$ is the highest power of $x_{q}$ dividing $x(a)$. Then by the construction of $D_q$ and
formulas (\ref{eqhe0}) and (\ref{chen1}), we know $x_{q}^{s_a}\nmid x(b)$. Note that $\triangle(b)\mid x(b)$.
Thus \begin{equation}\label{11111}x_{q}^{s_a}\nmid \triangle(b).\end{equation}
On the other hand, by statement (\ref{11110}), $x_q$ fulfills the conditions of Lemma \ref{sect3-0l}.
Thus,\begin{equation}\label{11112}x_{q}^{s_a} \mbox{ is the highest power of } x_{q} \mbox{ dividing }\triangle(a).\end{equation}
Therefore, $\triangle(a)\nmid \triangle(b)$ by (\ref{11111}).

Similarly, we can prove that $\triangle(b)\nmid \triangle(a)$.\\
$ $\\
B. Obviously, the map $g$ is meet-preserving.\\
$ $\\
C. The map $g$ is join-preserving.

Let $p, q\in P$. Obviously, $\mbox{supp}(p)\cup \mbox{supp}(q)\subseteq \mbox{supp}(p\vee q)$.
Now, let $$T_{p\vee q}=\mbox{supp}(p\vee q)\setminus (\mbox{supp}(p)\cup \mbox{supp}(q)).$$
Then $$g(p\vee q)=g(p)\vee g(q)\vee \mbox{lcm}\{\triangle (a_v): a_v\in T_{p\vee q}\}.$$ If $T_{p\vee q}=\emptyset$
then $g(p\vee q)=g(p)\vee g(q)\vee\mbox{lcm }\emptyset=g(p)\vee g(q)$.
Next, suppose that $T_{p\vee q}\neq \emptyset$. Then by Lemma \ref{sect3-1l},
$$\mbox{lcm}\{\triangle (a_v):a_v\in T_{p\vee q}\}\mid \mbox{lcm}\{\triangle (a_v): a_v\in \mbox{supp}(p)\cup \mbox{supp}(q)\}$$
since $\mbox{supp}(p)\cup \mbox{supp}(q) \in B_{p\vee q}$. Therefore, $g(p\vee q)=g(p)\vee g(q)$, i.e., the map $g$ is join-preserving.\\
$ $ \\
D. The map $g$ is surjective.

Assume that $p^{'}\in P^{'}$. Then $p^{'}= \mbox{lcm}\{\triangle(a_{i}): i\in I\}$ with $a_i\in \mbox{atoms}(P)$ for each $i\in I$. Let $b=\bigvee_{i\in I}a_{i}\in P$.
Then $\{a_i: i\in I\}\in B_b$. Thus, by Lemma \ref{sect3-1l}, $\triangle(a_{j})\mid \mbox{lcm}\{\triangle(a_{i}): i\in I\}$ for all $a_j\in \mbox{supp}(b)\setminus\{a_i: i\in I\}$.
Therefore, $$g(b)= \mbox{lcm}\{\triangle(a_{i}): a_{i}\in \mbox{supp}(b)\}=\mbox{lcm}\{\triangle(a_{i}): i\in I\}=p^{'},$$
which means that $g$ is surjective. \\
$ $\\
E. The map $g$ is injective.

Equivalently, we only need to prove that $a=b$ when $g(a)=g(b)$. For any $a, b\in P$, distinguishing two situations, we can have either $0\in \{a,b\}$ or $a, b\in P\setminus 0$. In the first case, we have $g(a)=g(b)=g(0)=1$. Obviously, $a=0=b$ by (\ref{11113}) and statement A.
In the second case, the proof will be completed by two parts.

(i) Suppose that $b\nleq a$. In this case, we easily see that
\begin{equation}\label{eq3.1}g(b)=\mbox{lcm}\{\triangle(a_{i}): a_{i}\in \mbox{supp}(b)\cap \mbox{supp}(a)\}\vee \mbox{lcm}\{\triangle(a_{j}): a_{j}\in \mbox{supp}(b)\setminus\mbox{supp}(a)\}.\end{equation}

From $b\nleq a$, $\mbox{supp}(b)\setminus\mbox{supp}(a)\neq \emptyset$. Now, let $a_r \in \mbox{supp}(b)\setminus\mbox{supp}(a)$.
Then $a_{r}\leq b$ but $a_{r}\nleq a$, which together with $a_r\in \mbox{atoms}(P)$ yields that
$a_{r}\| a$. Thus, by Lemma \ref{sect2-1l}, we have that $$\mbox{mi}(P)\cap\lceil a_{r}\rceil\nsubseteq \mbox{mi}(P)\cap\lceil a\rceil\mbox{ and }\mbox{mi}(P)\cap\lceil a\rceil\nsubseteq\mbox{mi}(P)\cap\lceil a_{r}\rceil,$$ and consequently
$$\mbox{mi}(P)\cap \lceil a\rceil^{c} \nsubseteq \mbox{mi}(P)\cap \lceil a_r\rceil^{c}
\mbox{ and }\mbox{mi}(P)\cap \lceil a_r\rceil^{c} \nsubseteq \mbox{mi}(P)\cap \lceil a\rceil^{c}.$$
Hence, there exists an element $q$ such that $q\in\mbox{mi}(P)\cap\lceil a_{r}\rceil^{c}$ but $q\notin\mbox{mi}(P)\cap\lceil a\rceil^{c}$.
Let $a_m\in \mbox{supp}(a)$. Then $a_m \leq a$ and which implies that $\lceil a_{m}\rceil^{c}\cap \mbox{mi}(P)\subseteq \lceil a\rceil^{c}\cap \mbox{mi}(P)$.
Thus, $q\notin \mbox{mi}(P)\cap\lceil a_{m}\rceil^{c}$. Therefore, \begin{equation}\label{xxxx}q\notin \mbox{mi}(P)\cap\lceil a_{i}\rceil^{c}\end{equation} for all $a_i\in \mbox{supp}(a)$.

By statement (\ref{11110}), there exists a variable $x_q$ in $m_q$ such that for all $r\in P$, $x_q\mid m_r$ implies that $q\nparallel r$.
Let $x_{q}^{s_{a_r}}$ be the highest power of $x_q$ dividing $x(a_r)$. Then similar to the proof of formula (\ref{11112}), we have $x_{q}^{s_{a_r}}\mid \triangle(a_{r})$. Thus
$$x_{q}^{s_{a_r}}\mid \mbox{lcm}\{\triangle(a_{j}): a_{j}\in \mbox{supp}(b)\setminus\mbox{supp}(a)\}$$
since $a_r \in \mbox{supp}(b)\setminus\mbox{supp}(a)$. Therefore, $x_{q}^{s_{a_r}}\mid g(b)$ by (\ref{eq3.1}).
Furthermore, similar to the proof of formula (\ref{11111}), from (\ref{xxxx}) we have that for all $a_i \in \mbox{supp}(a)$, $x_{q}^{s_{a_r}}\nmid \triangle(a_{i})$. Thus $x_{q}^{s_{a_r}}\nmid g(a)$. Consequently, $g(b)\nmid g(a)$, contrary to $g(a)=g(b)$.

Consequently, $b\leq a$.

(ii) Similar to the proof of (i), the condition $a\nleq b$ will deduce a contradiction.

With (i) and (ii) we know that $a=b$ if $g(a)=g(b)$ in the case that $a, b\in P\setminus 0$.

 Therefore, the map $g$ is injective.

From B, C, D and E, $g$ is an isomorphism from $P$ to $P^{'}$. Further, by (\ref{11113}), $\mathcal{M}$ is a weak coordinatization of $P$.
\endproof

The following two examples will illustrate Theorem \ref{sect3-1t}.
\begin{example}\label{sect3-2e1}
\emph{Let $P$ be a finite atomic lattice with a labeling as Fig.6. It is easy to see that the labeling of $P$ satisfies the conditions of Theorem \ref{sect3-1t} and does not satisfy the conditions of Lemma \ref{sect2-2l}.}
\par\noindent\vskip50pt
 \begin{minipage}{11pc}

\setlength{\unitlength}{0.75pt}\begin{picture}(100,100)

\put(65,20){\circle{4}}\put(55,15){\makebox(0,0)[l]{\footnotesize $a$}}

\put(125,20){\circle{4}}\put(115,15){\makebox(0,0)[l]{\footnotesize$e$}}

\put(185,20){\circle{4}}\put(200,15){\makebox(0,0)[r]{\footnotesize$m$}}

\put(65,80){\circle{4}}\put(45,85){\makebox(0,0)[l]{\footnotesize$ac$}}

\put(125,80){\circle{4}}\put(102,85){\makebox(0,0)[l]{\footnotesize$cm$}}

\put(185,80){\circle{4}}\put(195,85){\makebox(0,0)[r]{\footnotesize$e$}}

\put(125,140){\circle{4}}

  \put(65,21){\line(0,1){58}}
  \put(66,21){\line(1,1){58}}
  \put(124,21){\line(-1,1){58}}
  \put(126,21){\line(1,1){58}}
    \put(185,21){\line(0,1){58}}
  \put(184,21){\line(-1,1){58}}
  \put(66,81){\line(1,1){58}}
  \put(125,81){\line(0,1){58}}
  \put(184,81){\line(-1,1){58}}

 \put(10,-10){\emph{Fig.5. The lattice $P$ with labeling $\mathcal{M}$}}

 \put(310,20){\circle{4}}

\put(370,20){\circle{4}}

\put(430,20){\circle{4}}

\put(310,80){\circle{4}}

\put(370,80){\circle{4}}

\put(430,80){\circle{4}}

\put(370,140){\circle{4}}

  \put(310,21){\line(0,1){58}}
  \put(311,21){\line(1,1){58}}
  \put(369,21){\line(-1,1){58}}
  \put(371,21){\line(1,1){58}}
    \put(430,21){\line(0,1){58}}
  \put(429,21){\line(-1,1){58}}
  \put(311,81){\line(1,1){58}}
  \put(370,81){\line(0,1){58}}
  \put(429,81){\line(-1,1){58}}
 \put(310,-10){\emph{Fig.6. $LCM(I_{P,\mathcal{M}})$}}
 \end{picture}
 \end{minipage}
 \quad\\
\quad\\

\emph{One can clarify that $I_{P,\mathcal{M}}=(e^{2}m, acm^{2}, a^{2}ce)$, and $LCM(I_{P,\mathcal{M}})$ is isomorphic to $P$ (see Figs.5 and 6).
Moreover, one can check that $\mathcal{M}$ is a weak coordinatization and $I_{P,\mathcal{M}}= M_{P,\mathcal{M}}$.}
 \end{example}
\begin{example}\label{sect3-2e}
\emph{Let us consider the finite atomic lattice $P$ with a labeling as Fig.2 again. One can clarify that the labeling of $P$ satisfies the conditions
of Theorem \ref{sect3-1t} and does not satisfy the conditions of Lemma \ref{sect2-2l}. Moreover, the
labeling $\mathcal{M}$ is a weak coordinatization and $I_{P,\mathcal{M}}\neq M_{P,\mathcal{M}}$.}
 \end{example}

\begin{remark}\label{re 3.1}
\emph{From Theorem \ref{sect3-2t}, if the monomial ideal $M=M_{P,\mathcal{M}}$ with the labeling $\mathcal{M}$ satisfying the conditions of Lemma \ref{sect2-2l}, then $\mathcal{D}=\mathcal{M}$. On the other hand, by Lemma \ref{sect2-4l} we know that if the monomial ideal $M=I_{P, \mathcal{M}}$ with the labeling $\mathcal{M}$ satisfies the conditions of Theorem \ref{sect3-1t} and does not satisfy the conditions of Lemma \ref{sect2-2l}, then $M$ must induce a new labeling $\mathcal{D}$ which is different from $\mathcal{M}$ and $D_{P,\mathcal{D}}=I_{P, \mathcal{M}}=M$.}
\end{remark}

\section{Finite super-atomic lattices}
Let $\mathcal{L}(n)$ be the set of all finite atomic lattices with $n$ ordered atoms. $\mathcal{L}(n)$ has a partial order where $Q\leq P$ if and only if there exists a join-preserving map which is a bijection on atoms from $P$ to $Q$ (note that such a map will also be surjective)(see \cite{Mapes13}). In this section, we shall discuss the structure of lattice $\mathcal{L}(n)$. We shall first define a finite super-atomic lattice, and then give an algorithm to find out all the finite super-atomic lattices in $\mathcal{L}(n)$.

\begin{definition}\label{sect4-1d}
 \emph{A finite atomic lattice $P$ is called super-atomic if it satisfies that for each $p\in (P\setminus\mbox{atoms}(P))\setminus 0$,
 there exists $T_0=\{a_1, a_2\}\in B_p$ such that $T_0\subseteq T$ for any $T\in B_p$.}
\end{definition}

For example, the finite atomic lattice $P$ shown as Fig.7 is super-atomic.
\par\noindent\vskip80pt
 \begin{minipage}{11pc}
\setlength{\unitlength}{0.75pt}\begin{picture}(300,120)

\put(100,40){\circle{4}}

\put(200,40){\circle{4}}

\put(300,40){\circle{4}}

\put(400,40){\circle{4}}

\put(150,90){\circle{4}}
\put(250,90){\circle{4}}
\put(350,90){\circle{4}}
\put(200,140){\circle{4}}
\put(250,190){\circle{4}}
\put(300,140){\circle{4}}
  \put(101,41){\line(1,1){47}}
  \put(201,41){\line(1,1){47}}
  \put(199,41){\line(-1,1){47}}
  \put(299,41){\line(-1,1){47}}
  \put(399,41){\line(-1,1){47}}
  \put(201.1,40.8){\line(3,1){147.2}}
  \put(151,91){\line(1,1){47}}
 \put(249,91){\line(-1,1){47}}
  \put(349,91){\line(-1,1){47.3}}
  \put(201,141){\line(1,1){47}}
  \put(299,141){\line(-1,1){47}}
  \put(251,91){\line(1,1){47}}
 \put(130,0){Fig.7. A finite super-atomic lattice}
 \end{picture}
 \end{minipage}

\begin{theorem}\label{sect5-1t}
A lattice $P$ is super-atomic if and only if for each $p\in (P\setminus\emph{atoms}(P))\setminus 0$,
there exits $\{ a_1,  a_2\}\in B_p$ such that $\emph{supp(p)}\setminus a_1\in\mathcal{S}_{P}$
and $\emph{supp}(p)\setminus a_2\in\mathcal{S}_{P}$.
\end{theorem}
\proof
Suppose that $P$ is super-atomic. Then there exits $\{ a_1,  a_2\}\in B_p$ for each $p\in (P\setminus\mbox{atoms}(P))\setminus 0$.
Now, assume that $\mbox{supp}(p)\setminus a_1\notin \mathcal{S}_{P}$.
Then, $$\mbox{supp}(p)\supseteq\overline{\bigvee}_{a\in \mbox{supp}(p)\setminus a_1} \mbox{supp}(a)=
\overline{\bigvee}_{a\in \mbox{supp}(p)\setminus a_1} \{a\}\supsetneq \mbox{supp}(p)\setminus a_1,$$
in which $\overline{\bigvee}$ is the join of $(\mathcal{S}_{P}, \subseteq)$.
Thus \begin{equation}\label{111001}\overline{\bigvee}_{a\in \mbox{supp}(p)\setminus a_1} \mbox{supp}(a)= \mbox{supp}(p).\end{equation}
By the definition of $\mathcal{S}_{P}$, $(\mathcal{S}_{P}, \subseteq)$ is the same as lattice $P$. Thus
\begin{equation}\label{1112}\mbox{supp}(q) \mbox{ corresponds to } q\mbox{ for each }q\in P,\end{equation}
and
\begin{equation}\label{11102}\mbox{ for any }S\in \mathcal{S}_{P}, \mbox{ there exists } q\in P \mbox{ such that } S=\mbox{supp}(q).\end{equation}
Therefore, by formulas (\ref{111001}) and (\ref{1112}), we have $\bigvee_{a\in \mbox{supp}(p)\setminus a_1} a= p$,
which means that $\mbox{supp}(p)\setminus a_1 \in B_p$. As $P$ is super-atomic,
there exists $T_0=\{b_1, b_2\}\in B_P$ such that $T_0\subseteq\mbox{supp}(p)\setminus a_1\bigcap \{a_1, a_2\}$, a contradiction.
So that $\mbox{supp}(p)\setminus a_1\in \mathcal{S}_{P}$.
Similarly, we can prove that $\mbox{supp}(p)\setminus a_2\in \mathcal{S}_{P}$.

Conversely, let $p\in (P\setminus\mbox{atoms}(P))\setminus 0$.
Then by the hypothesis, there exists $\{ a_1,  a_2\}\in B_p$ such that $\mbox{supp}(p)\setminus a_1\in\mathcal{S}_{P}$
and $\mbox{supp}(p)\setminus a_2\in\mathcal{S}_{P}$. Note that $T\subseteq \mbox{supp}(p)$ for all $T\in B_p$.
Next, we prove that $\{a_1, a_2\}\subseteq T $ for all $T\in B_p$. If there exists a $T\in B_p$ such that $\{a_1, a_2\}\nsubseteq T$, then either
$\overline{\bigvee}_{a\in T}\{a\}=\overline{\bigvee}_{a\in T}\mbox{supp}(a)\subseteq \mbox{supp}(p)\setminus a_1\in \mathcal{S}_{P}$
or $\overline{\bigvee}_{a\in T}\{a\}=\overline{\bigvee}_{a\in T}\mbox{supp}(a)\subseteq \mbox{supp}(p)\setminus a_2 \in \mathcal{S}_{P}$.
By (\ref{1112}) and (\ref{11102}), in any case we have that $\bigvee_{a\in T}a<p$, contrary to $T\in B_p$.
Hence, \begin{equation}\label{0005}\{a_1, a_2\}\subseteq T \mbox{ for all } T\in B_p.\end{equation}

Therefore, by Definition \ref{sect4-1d} and (\ref{0005}), $P$ is a finite super-atomic lattice.
\endproof

By Definition \ref{sect4-1d} and Theorem \ref{sect5-1t}, we have the lemma as below, obviously.

\begin{lemma}\label{sect4-1p}
Let $P$ be a super-atomic lattice in $\mathcal{L}(n)$ with $\emph{atoms}(P)=\{1, 2, \cdots, n\}$ and $n\geq 2$. Then $(\mathcal{S}_{P}, \subseteq)$ satisfies the following statements:\\
\emph{(D1)} $\{\emptyset, \{1\}, \cdots, \{n\}, \{1, \cdots, n\}\}\subseteq \mathcal{S}_P$. \\
\emph{(D2)} If $S\in \mathcal{S}_P\setminus\{\emptyset, \{1\}, \cdots, \{n\}\}$, then
there exist two different atoms $\{i\},\{j\}\in \mathcal{S}_P$ such that $S=\{i\}\vee \{j\} $ and $S\setminus k \in \mathcal{S}_P$ for any $k\in \{i,j\}$.\\
\emph{(D3)} Let $S_1, S_2 \in \mathcal{S}_P$. If $S_1=\{u\}\vee \{v\}$, $S_2=\{k\}\vee \{h\}$ and $S_1\|S_2$, then $\{u,v\}\nsubseteq S_2$ and $\{k,h\}\nsubseteq S_1$.
\end{lemma}

In what follows, we shall suggest an algorithm to construct all finite super-atomic lattices in $\mathcal{L}(n)$ with $n\geq 2$.
\begin{algorithm}\label{00013}\quad\\
\textbf{Input:} $X=\{1,2,\cdots, n\}$.\\
\textbf{Output:} $\mathcal{S}^*$.\\
Step 1. Take $\mathcal{S}_0=\{\emptyset\}, \mathcal{S}_1=\{\{1\}, \cdots, \{n\}\}, \mathcal{S}_n=\{X\}$,
$\mathcal{S}^*=:\mathcal{S}_0 \cup \mathcal{S}_1 \cup \mathcal{S}_n$ and $k:=0$.\\
Step 2. If $n-k=2$, then go to Step 7.\\
Step 3. For any $S\in \mathcal{S}_{n-k}$, take $\delta(S)=\{i_S, j_S\}\subseteq S$ satisfying $\delta(S)\nsubseteq T$ for all $T\in \mathcal{S}_{n-k}\setminus S$.\\
Step 4. $\mathcal{S}_{n-k-1}=\bigcup_{S\in \mathcal{S}_{n-k}}\{S\setminus i_S, S\setminus j_S\}$.\\
Step 5. $k:=k+1$.\\
Step 6. $\mathcal{S}^*:=\mathcal{S}^* \cup \mathcal{S}_{n-k}$, and go to Step 2.\\
Step 7. Stop.
\end{algorithm}

\begin{theorem}\label{t-01}
Every output $(\mathcal{S}^{*},\subseteq )$ in Algorithm \ref{00013} is a finite super-atomic lattice in $\mathcal{L}(n)$.
Further, every finite super-atomic lattice in $\mathcal{L}(n)$ can be constructed by Algorithm \ref{00013}.

\end{theorem}
\proof  Throughout the proof, let $\bigvee \delta(S)=\{i_S\}\vee \{j_S\}$ for any $S\in \mathcal{S}^*\setminus (\mathcal{S}_0 \bigcup \mathcal{S}_1)$.
First, we shall prove that every output $(\mathcal{S}^{*},\subseteq )$ in Algorithm \ref{00013} is a finite super-atomic lattice by four steps as below.\\
\quad\\
(B1). Obviously, $(\mathcal{S}^{*},\subseteq )$ has a minimum element $\emptyset$
and a maximum element $\{1,\cdots, n\}$.\\
\quad\\
(B2). If $S\in \mathcal{S}^{*}\setminus(\mathcal{S}_1\bigcup \mathcal{S}_0)$ then $S=\bigvee \delta(S)$.\\
\quad

Observe that there exists $t\in \{2, \cdots, n\}$ such that $S\in\mathcal{S}_t$
and \begin{equation}\label{00014}\delta(S)\nsubseteq T\end{equation} for all $T\in \mathcal{S}_{t}\setminus S$ by Algorithm \ref{00013}.
Set \begin{equation}\label{00015}\mathcal{D}= \{D\in \mathcal{S}^{*}: \delta(S) \subseteq D\} \mbox{ and } \mathcal{D}_{\ast}=
\{D: D \mbox{ is a minimal element of  }\mathcal{D}\}.\end{equation}
Let $D\in \mathcal{D}_{\ast}$. We claim that $D\notin \mathcal{S}_{u}$ for any integer $u$ with $0\leq u< t$.
Indeed, if $D\in \mathcal{S}_{u}$, then there exists $G\in \mathcal{S}_{t}$ such that $D\subsetneq G$
by Algorithm \ref{00013}. Thus $\delta(S) \subseteq G$, which together with (\ref{00014}) yields that $G=S$.
Therefore, $D\subseteq S\setminus i_S$ or $D\subseteq S\setminus j_S$ by Algorithm \ref{00013},
contrary to $\delta(S) \subseteq D$.

Below, assume that $D\in \mathcal{S}_v$ with $n\geq v\geq t$. Now, we shall prove $v=t$.
Suppose that $n\geq v> t$. By Algorithm \ref{00013}, there exists $R\in \mathcal{S}_v$ such that $R\supsetneq S$. There are two cases.

Case (1). If $D=R$ then $D \supsetneq S$, contrary to $D\in \mathcal{D}_{\ast}$ since $S\in \mathcal{D}$.

Case (2). Let $D\neq R$. We first claim that $\delta(D)=\delta(S)$. Otherwise, either $\delta(S)\subseteq D\setminus i_D\subsetneq D$
or $\delta(S)\subseteq D\setminus j_D\subsetneq D$, contrary to $D\in \mathcal{D}_{\ast}$.
Hence, $\delta(D)= \delta(S) \subseteq S\subsetneq R$, contrary to $\delta(D) \nsubseteq R$ since $R\neq D$
and both $R$ and $D$ in $\mathcal{S}_v$ (see (\ref{00014})).

Cases (1) and (2) imply that $v=t$. Therefore, $D=S$ by formulas (\ref{00014}) and (\ref{00015}),
which means that $\mathcal{D}_{\ast}$ contains exactly one element $S$ and $S=\bigvee\delta(S)$.\\
\quad\\
(B3). If $S_1, S_2 \in \mathcal{S}^*$ then $S_1 \vee S_2$ exists in $\mathcal{S}^{*}$.\\
\quad

Obviously, if $S_1\nparallel S_2$ then $S_1\vee S_2=S_1$ or $S_1\vee S_2=S_2$.

Next, suppose that $S_1 \| S_2$.
Observe that $S_1$ and $S_2$ are not in $\mathcal{S}_0$. There are three cases.

Case (i). If $S_1=\{i\}, S_2=\{j\}$ and $i\neq j$, then $S_1 \vee S_2$ exists in $\mathcal{S}^{*}$.

In this case, set
$$M= \{S\in \mathcal{S}^{*}: \{i,j\} \subseteq S\} \mbox{ and } M_{\ast}= \{S: S \mbox{ is a minimal element of  }M\}.$$
Note that $M \neq \emptyset$. Hence $M_{\ast}\neq \emptyset$.
Assume that $S\in M_{\ast}$. Then $S\in \mathcal{S}^{*}\setminus(\mathcal{S}_1\bigcup \mathcal{S}_0)$. Thus by (B2),
$S= \bigvee \delta(S)$. If $\{i,j\}\neq \delta(S)$ then $\{i,j\} \subseteq S\setminus i_S \in \mathcal{S}^*$ or $\{i,j\} \subseteq S\setminus j_S\in \mathcal{S}^*$
by Algorithm \ref{00013}, contrary to the fact that $S\in M_*$. Therefore, $\{i,j\}= \delta(S)$, which means that $S_1\vee S_2= S\in \mathcal{S}^*$.

Case (ii). If $S_1=\{i\}$ and $S_2\in \mathcal{S}^{*}\setminus(\mathcal{S}_1\bigcup \mathcal{S}_0)$ with $i\notin S_2$, then $S_1 \vee S_2$ exists in $\mathcal{S}^{*}$.

Indeed, by (B2), $S_2= \bigvee \delta(S_2)$.
Suppose that $S_1 \vee S_2$ dose not exist in $\mathcal{S}^*$. Then $\mathcal{S}^*$ contains two different
minimal elements containing $S_1 \bigcup S_2$, say $S_a, S_b$. Clearly, $S_a\| S_b$.

We claim that \begin{equation}\label{000022}\delta(S_a) \subseteq \{i, i_{S_2}, j_{S_2}\}\mbox{ and }\delta(S_a)\neq \delta(S_2).\end{equation}
Suppose that $\delta(S_a) \nsubseteq \{i, i_{S_2}, j_{S_2}\}$. By Algorithm \ref{00013}, $\{i, i_{S_2}, j_{S_2}\}\subseteq S_a\setminus i_{S_a}\in \mathcal{S}^*$
or $\{i, i_{S_2}, j_{S_2}\}\subseteq S_a\setminus j_{S_a}\in \mathcal{S}^* $.
From $S_2= \bigvee \delta(S_2)$, if $\{i, i_{S_2}, j_{S_2}\}\subseteq S_a\setminus i_{S_a}$ then $S_1 \bigcup S_2 \subseteq S_a\setminus i_{S_a} \subsetneq S_a$, a contradiction.
Similarly, we can prove that $\{i, i_{S_2}, j_{S_2}\}\subseteq S_a\setminus j_{S_a}\in \mathcal{S}^* $ will deduce a contradiction. Therefore, $\delta(S_a) \subseteq \{i, i_{S_2}, j_{S_2}\}$. Now assume that $\delta(S_a)=\delta(S_2)$. Then $S_a = \bigvee\delta(S_a)=\bigvee\delta(S_2)=S_2$ which implies $i\in S_2$, a contradiction.

Arguing as formula (\ref{000022}), we have \begin{equation}\label{000023}\delta(S_b) \subseteq \{i, i_{S_2}, j_{S_2}\}\mbox{ and }\delta(S_b)\neq \delta(S_2).\end{equation}

Formulas (\ref{000022}) and (\ref{000023}) imply that both $\delta(S_a)$ and $\delta(S_b)$ equal to $\{i, i_{S_2}\}$ or $\{i, j_{S_2}\}$.
We claim that \begin{equation}\label{yyy}\delta(S_a)\neq \delta(S_b).\end{equation} Indeed if $\delta(S_a)=\delta(S_b)$ then $S_a=\bigvee\delta(S_a)=\bigvee\delta(S_b)=S_b$, contrary to $S_a\|S_b$.
Thus, if $\delta(S_a)=\{i, i_{S_2}\}$ then $\delta(S_b)=\{i, j_{S_2}\}$.
Clearly, $\{i, j_{S_2}\}\subseteq S_a\setminus i_{S_2} \in \mathcal{S}^*$. Thus $S_b=\bigvee \delta(S_b)= \{i\}\vee \{j_{S_2}\} \subseteq S_a\setminus i_{S_2} \subsetneq S_a$,
contrary to $S_a\| S_b$. Similarly, we can prove that $\delta(S_a) =\{i, j_{S_2}\}$ will deduce a contradiction.
Therefore, $S_1 \vee S_2$ exists in $\mathcal{S}^{*}$.

Case (iii). If $S_1, S_2\in \mathcal{S}^{*}\setminus(\mathcal{S}_1\bigcup \mathcal{S}_0)$ and $S_1\|S_2$, then $S_1 \vee S_2$ exists in $\mathcal{S}^{*}$.

First, if $\delta(S_1)\subseteq S_2$
then $\bigvee \delta(S_1)=S_1 \subseteq S_2$, a contradiction. Thus $\delta(S_1)\nsubseteq S_2$. Similarly, we can prove $\delta(S_2)\nsubseteq S_1$.

Then assume that $S_1 \vee S_2$ does not exist in $\mathcal{S}^*$. Then $\mathcal{S}^*$ contains two different
minimal elements containing $S_1 \bigcup S_2$, say $C_1, C_2$. Clearly, $C_1 \| C_2$.
Similar to the proof of formula (\ref{000022}) in Case (ii), we can prove that \begin{equation} \label{00018}\delta(C_1)\subseteq \delta(S_1) \cup\delta(S_2), \delta(C_1)\neq\delta(S_1)\mbox{ and }\delta(C_1)\neq\delta(S_2).\end{equation}

Using (\ref{00018}), we know that $\delta(C_1)$ equals to one of four sets $\{i_{S_1}, i_{S_2}\}, \{i_{S_1}, j_{S_2}\}, \{j_{S_1}, i_{S_2}\}$
and $\{j_{S_1}, j_{S_2}\}$. Similarly, we can prove that $\delta(C_2)$ also equals to one of four sets $\{i_{S_1}, i_{S_2}\}$, $\{i_{S_1}, j_{S_2}\}$, $\{j_{S_1}, i_{S_2}\}$
and $\{j_{S_1}, j_{S_2}\}$.  Similar to the proof of formula (\ref{yyy}) in Case (ii), we can prove $\delta(C_1)\neq \delta(C_2)$.
Now, suppose that $\delta(C_1)=\{i_{S_1}, i_{S_2}\}$. Then $C_1\setminus i_{S_1}, C_1\setminus i_{S_2} \in \mathcal{S}^{*}$ by Algorithm \ref{00013}.
If $\delta(C_2)=\{i_{S_1}, j_{S_2}\}$ then $C_2=\{i_{S_1}\}\vee \{j_{S_2}\} \subseteq C_1\setminus i_{S_2}$, contrary to $C_1 \| C_2$.
Similarly, we can prove that all the other cases will deduce a contradiction. Therefore, $S_1 \vee S_2$ exists in $\mathcal{S}^{*}$.\\
\quad\\
(B4). $(\mathcal{S}^{*},\subseteq )$ is super-atomic.\\

By (B1), (B2), (B3) and  Definitions \ref{sect2-2d} and \ref{sect2-3d}, $(\mathcal{S}^{*}, \subseteq)$ is a finite atomic lattice.
Next, we shall prove that $(\mathcal{S}^{*}, \subseteq)$ is super-atomic.

Suppose $S\in \mathcal{S}^{*}\setminus(\mathcal{S}_1\bigcup \mathcal{S}_0)$ and $T\in B_S$. Note that $\bigvee T=S$.
If $\{i_{S}\} \notin T$ then $\bigcup T\subseteq S\setminus i_{S} \in \mathcal{S}^*$, which implies $\bigvee T\subseteq S\setminus i_{S}$, contrary to $\bigvee T=S$.
Thus $\{i_{S}\}\in T$. Similarly, we have $\{j_{S}\}\in T$. Hence $\{\{i_S\}, \{j_S\}\}\subseteq T$. Again by (B2), $\{i_{S}\}\vee \{j_{S}\} =\bigvee\delta(S)=S$, this means that
$\{\{i_S\}, \{j_S\}\}\in B_S$. Thus, by Definition \ref{sect4-1d}, the lattice $(\mathcal{S}^*, \subseteq)$ is super-atomic.

$ $

We finally prove that every super-atomic lattice in $\mathcal{L}(n)$ can be constructed by Algorithm \ref{00013}.

Let $(\mathcal{S}, \subseteq)$ be a super-atomic lattice in $\mathcal{L}(n)$. For each $0\leq i\leq n$, define $\mathcal{T}_i=\{S\in \mathcal{S}: |S|=i\}$. Then $\mathcal{S}=\mathcal{T}_0 \bigcup \mathcal{T}_1\bigcup \cdots \bigcup \mathcal{T}_n$. In what follows, we prove that there is an output $\mathcal{S}^{*}$ by Algorithm \ref{00013}
such that $\mathcal{S}^{*}=\mathcal{S}$.

In fact, from Algorithm \ref{00013}, we know that $\mathcal{S}^{*}=\mathcal{S}_0 \bigcup \mathcal{S}_1\bigcup \cdots \bigcup \mathcal{S}_n$. Therefore, in order to construct $\mathcal{S}^{*}$ by Algorithm \ref{00013} such that $\mathcal{S}^{*}=\mathcal{S}$, we just need to construct $\mathcal{S}_i$ such that $\mathcal{T}_i=\mathcal{S}_i$ for all $0\leq i\leq n$.

First, by Algorithm \ref{00013} and (D1) in Lemma \ref{sect4-1p}, we have \begin{equation}\label{22221}\mathcal{T}_i=\mathcal{S}_i \mbox{ for all }i\in \{0, 1, n\}.\end{equation}

Then by (D2), there exist $\{i^S\}, \{j^S\}\in\mathcal{T}_1$ such that $\{i^S\}\vee \{j^S\}=S$ for each $S\in \mathcal{T}_n$ with $n\geq 2$.
As $\mathcal{T}_n=\mathcal{S}_n$, by (\ref{22221}), we can take $\delta(S)=\{i^S, j^S\}$ in Step 3 of Algorithm \ref{00013} for all $S\in \mathcal{S}_n$. Thus $\mathcal{T}_{n-1}\supseteq\mathcal{S}_{n-1}$ by
Step 4 and (D2).

We claim that $\mathcal{T}_{n-1}=\mathcal{S}_{n-1}$.
Otherwise, there exists $W\in \mathcal{T}_{n-1}$ such that $W\notin \mathcal{S}_{n-1}$. Let $K\in \mathcal{S}$ with $K\succ W$.
Then by (D2), there exist $\{i^K\}, \{j^K\}\in\mathcal{T}_1$ such that $\{i^K\}\vee \{j^K\}=K$. If $\{i^K, j^K\}\subseteq W$
then $K=\{i^K\}\vee \{j^K\}\subseteq W\prec K$, a contradiction. Thus $\{i^K, j^K\}\nsubseteq W$. It follows (D2) that $W\subseteq K\setminus i^K \prec K$
or $W\subseteq K\setminus j^K\prec K$, which means that $W=K\setminus i^K$ or $W=K\setminus j^K$.
Therefore, $K\in \mathcal{T}_n$, which together with $\mathcal{T}_n=\mathcal{S}_n$ yields that
$W\in \mathcal{S}_{n-1}$ since $\delta(K)=\{i^K, j^K\}$, a contradiction.

Similarly, we can construct $\mathcal{T}_{h}=\mathcal{S}_{h}$ by taking $\delta(T)=\{i^T, j^T\}$ for any $T\in \mathcal{T}_{h+1}$,
in which  $\{i^T\}, \{j^T\}\in\mathcal{T}_1$ and $\{i^T\}\vee \{j^T\}=T$ for all $2\leq h\leq n-2$.

Consequently, $\mathcal{T}_i=\mathcal{S}_i$ for all $0\leq i\leq n$.
\endproof

The following example will illustrate Algorithm \ref{00013}.
\begin{example}\label{sect5-11e}
\emph{Let $n=3$. Then by Algorithm \ref{00013} we have three super-atomic lattices in $\mathcal{L}(n)$ as follows.
$$\mathcal{Q}_{1}=\{\emptyset, \{1\}, \{2\}, \{3\}, \{1,2,3\}, \{1,2\}, \{1,3\}\},$$
$$\mathcal{Q}_{2}=\{\emptyset, \{1\}, \{2\}, \{3\}, \{1,2,3\}, \{1,2\}, \{2,3\}\},$$
$$\mathcal{Q}_{3}=\{\emptyset, \{1\}, \{2\}, \{3\}, \{1,2,3\}, \{1,3\}, \{2,3\}\}.$$
On can check that $(\mathcal{Q}_{1}, \subseteq)$, $(\mathcal{Q}_{2}, \subseteq)$ and $(\mathcal{Q}_{3}, \subseteq)$
are all the super-atomic lattices in $\mathcal{L}(n)$.}
\end{example}

\section{Specific labelings}
In \cite{Mapes09}, there are three specific coordinatizations, i.e., Minimal Squarefree, Minimal Depolarized and Greedy, one can see that all of them are based on the labeling described as in Lemma \ref{sect2-2l}. In this section, we shall give a kind of labelings on a lattice $P$ which does not satisfy the conditions of Lemma \ref{sect2-2l}, and show the conditions that our labeling is either a coordinatization or a weak coordinatization.

Let $P\in \mathcal{L}(n)$ with $\mbox{atoms}(P)=\{a_{1}, a_{2},\cdots, a_{n}\}$. We define a labeling
$\mathcal{C}$ of $P$ as that $\mathcal{C}=\{m_p: p\in P\setminus 0\}$ where
\begin{equation}\label{eq5.2}m_p=\prod_{a_{i}\in \mbox{supp}(p)}a_{i}\end{equation}
in which every $a_i$ means both atom in $P$ and variable in labeling
$\mathcal{C}$.

In what follows, let $[a, b]=\{p\in P: a\leq p\leq b\}$ and $N([a, b])= |[a, b]|$ for the purposes of convenience.

\begin{theorem}\label{sect4-1t}
Let $P\in \mathcal{L}(n)$. For each $p\in (P\setminus\emph{atoms}(P))\setminus0$, if there exist $a_{i}, a_{j}\in \emph{supp}(p)$ such that $p=a_{i}\vee a_{j}$ and $N([a_{r}\vee a_{k}, 1])< N([p, 1])$ for a fixed number $r\in \{i,j\}$ and all $a_{k}\in \emph{atoms}(P)\setminus\emph{supp}(p)$, then the labeling $\mathcal{C}$ of $P$ as defined by (\ref{eq5.2}) is a weak coordinatization.
\end{theorem}
\proof
For $b\in P$, define $g: P\longrightarrow LCM(I_{P, \mathcal{C}})$ to be a map such that
$$g(b)= \mbox{lcm}\{\triangle(u): u\in \mbox{supp}(b)\}.$$
The main part is to show that $g$ is an isomorphism of lattices. Similar to B, C and D in the proof of Theorem \ref{sect3-1t}, we can prove that the map $g$ is meet-preserving , join-preserving and surjection. Thus, we only need to show that $g$ is injective. The proof will be split into two parts.

\quad\\
($\ast$) Let $a_u, a_v \in \mbox{atoms}(P)$. Then $a_u\mid \triangle(a_v)$ if and only if $a_u\neq a_v$.

Suppose that $a_u\mid \triangle(a_v)$. From formula (\ref{eq5.2}), $a_u\mid m_{p}$ if and only if $p\geq a_u$. Thus $a_u\nmid x(a_u)$ by (\ref{eqhe0}). This means that $a_u\nmid \triangle(a_u)$ since $\triangle(a_u) \mid x(a_u)$. Therefore, $a_u\neq a_v$.

Conversely, assume that $a_w \in \mbox{atoms}(P)\setminus a_u$. Then $a_u\in \lceil a_w\rceil^{c}$. Thus $a_u\mid x(a_w)$ by equations (\ref{eqhe0}) and (\ref{eq5.2}).
On the other hand, let $F\in \bigcup_{p\geq a_v}B_p$. Then $\bigvee F\geq a_v$. So that $a_u\neq\bigvee F$ since $a_u\neq a_v$. Thus there exists $a_z\in F$ such that $a_z\neq a_u$.
Hence $a_u\mid \mbox{lcm}\{x(a_j): a_j\in F\}$. This together with equation (\ref{equ1.2}) implies that $a_u\mid \triangle(a_v)$.

\quad\\
($\ast\ast$) The map $g$ is injective.

Clearly, if $0\in \{a, b\}$ and $g(a)=g(b)$ then $g(a)=g(b)=g(0)=1$, which implies that $a=0=b$.
Next, let $a, b\in P\setminus 0$ and $g(a)=g(b)$. Now we shall prove $a=b$.

Suppose that $b\nleq a$. Then we have either $a\in \mbox{atoms}(P)$ or $a\in (P\setminus\mbox{atoms}(P))\setminus0$.
In the first case, we have $\mbox{supp}(b)\setminus\mbox{supp}(a)\neq\emptyset$. Thus there exists $c\in \mbox{supp}(b)\setminus\mbox{supp}(a)$.
By statement ($\ast$), $a \nmid \triangle(a)$ and $a \mid \triangle(c)$. Therefore, $a\mid g(b)$ and $a\nmid g(a)$, a contradiction.

In the second case, let $a_k\in \mbox{atoms}(P)\setminus\mbox{supp}(a)$.
Then by the hypothesis of the theorem, there exist two elements $a_{i}, a_{j}\in \mbox{supp}(a)$ such that
$a=a_{i}\vee a_{j}$ and $N([a_{j}\vee a_{k}, 1])< N([a, 1])$ (set $r=j$).
For convenience, let $a_{j}^{n_{y}}$ be the highest power of $a_{j}$ dividing $x(a_{y})$ for each $a_y\in \mbox{atoms}(P)$. Clearly, by (\ref{eqhe0})
$$x(a_k)=\prod_{q\in\lceil a_k\rceil^{c}}m_{q}=\prod_{q_1\in\lceil a_k\rceil^{c}\cap\lceil a_j\rceil} m_{q_1}\ast\prod_{q_2\in\lceil a_k\rceil^{c}\cap\lceil a_j\rceil^{c}} m_{q_2}.$$
Thus by (\ref{eq5.2}), $n_k=|\lceil a_k\rceil^{c}\cap\lceil a_j\rceil|$.
On the other hand, $\lceil a_k\rceil^{c}\cap\lceil a_j\rceil=[a_j, 1]-[a_j\vee a_k, 1]$. So that $n_k=N([a_j, 1])-N([a_j\vee a_k, 1])$.
Similarly, $n_i=N([a_j, 1])-N([a_j\vee a_i, 1])$. Therefore,
\begin{equation}\label{eq5.1}n_k-n_i=N([a_j\vee a_i, 1])-N([a_j\vee a_k, 1])=N([a, 1])-N([a_j\vee a_k, 1])\geq 1.\end{equation}

Let $r\geq a_k$. Suppose $T\in B_r$. We claim that there exists $a_t\in T$ such that $a_t\in \mbox{atoms}(P)\setminus\mbox{supp}(a)$. Otherwise, $T\subseteq \mbox{supp}(a)$,
which means that $a_k\leq r=\bigvee T \leq \bigvee\mbox{supp}(a)= a$, contrary to $a_k\notin \mbox{supp}(a)$. Hence, $n_t-n_i\geq 1$ by (\ref{eq5.1}). Thus \begin{equation}\label{000024}a_{j}^{n_{i}+1}\mid \mbox{lcm}\{x(a_w): a_w\in T\}.\end{equation}
Below, let $a_{j}^{m_{y}}$ be the highest power of $a_{j}$ dividing $\triangle(a_{y})$ for each $a_y\in \mbox{atoms}(P)$. Thus $m_k\geq n_i+1$ by formulas (\ref{000024}) and (\ref{equ1.2}). Clearly, $m_i \leq n_i$ since $\triangle(a_{i})\mid x(a_i)$. Therefore \begin{equation}\label{000025}m_k> n_i\geq m_i.\end{equation}

Clearly, there exists $a_s\in\mbox{supp}(b)\setminus\mbox{supp}(a)\subseteq \mbox{atoms}(P)\setminus\mbox{supp}(a)$
such that $a_s\vee a_e=b$ for some $a_e\in \mbox{supp}(b)$. This follows that $g(b)=\mbox{lcm}\{\triangle(a_s),\triangle(a_e)\}$
since $g$ is join-preserving. Now, let $a_{j}^{m}$ be the highest power of $a_{j}$ dividing $g(b)$. Then $m\geq m_s$.
Using formula (\ref{000025}), $m_s> n_i\geq m_i$. Thus $m\geq m_s > m_i$.

On the other hand, $g(a)=\mbox{lcm}\{\triangle(a_{i}),\triangle(a_{j})\}$. By statement ($\ast$), we have $a_j\nmid \triangle(a_{j})$. Thus $a_{j}^{m_{i}}$ is
the highest power of $a_{j}$ dividing $g(a)$. As $m\geq m_s > m_i$, we finally have that $g(b)\nmid g(a)$, a contradiction.

Therefore, the assumption of $b\nleq a$ will deduce a contradiction. Consequently, $b\leq a$.

Similarly, we can prove that $a\leq b$, it follows from $b\leq a$ that $a=b$ finally.
\endproof

\begin{remark}\label{sect4-1r}
\emph{The labeling $\mathcal{C}$ as defined by (\ref{eq5.2}) needs not to satisfy the condition (C2) generally. For example, consider the lattice shown as Fig.8.}
\end{remark}
\par\noindent\vskip60pt
 \begin{minipage}{11pc}
\setlength{\unitlength}{0.75pt}\begin{picture}(300,160)

\put(100,40){\circle{4}}\put(80,28){\makebox(0,0)[l]{\footnotesize $a$}}

\put(200,40){\circle{4}}\put(180,28){\makebox(0,0)[l]{\footnotesize$b$}}

\put(300,40){\circle{4}}\put(320,28){\makebox(0,0)[r]{\footnotesize$c$}}

\put(400,40){\circle{4}}\put(420,28){\makebox(0,0)[r]{\footnotesize$d$}}

\put(150,90){\circle{4}}\put(130,100){\makebox(0,0)[l]{\footnotesize$ab$}}
\put(250,90){\circle{4}}\put(270,100){\makebox(0,0)[r]{\footnotesize$bc$}}
\put(350,90){\circle{4}}\put(370,100){\makebox(0,0)[r]{\footnotesize$ad$}}
\put(200,140){\circle{4}}\put(180,150){\makebox(0,0)[l]{\footnotesize$abc$}}
\put(250,190){\circle{4}}

  \put(101,41){\line(1,1){47}}
  \put(201,41){\line(1,1){47}}
  \put(199,41){\line(-1,1){47}}
  \put(299,41){\line(-1,1){47}}

  \put(399,41){\line(-1,1){47}}
  \put(101.3,40.5){\line(5,1){248}}
  \put(151,91){\line(1,1){47}}

  \put(249,91){\line(-1,1){47}}
  \put(349,91){\line(-1,1){98}}
  \put(201,141){\line(1,1){47}}

 \put(115,5){Fig.8. The lattice $P$ with a labeling $\mathcal{C}$}
 \end{picture}
 \end{minipage}
 \\
Clearly, the lattice $P$ satisfies the conditions of Theorem \ref{sect4-1t}, and its labeling $\mathcal{C}$ yields
that $I_{P,\mathcal{C}}= \{b^{2}c^{2}d, a^{2}cd^{2}, a^{3}b^{2}d^{2}, a^{3}b^{4}c^{3}\}$. On can check that $LCM(I_{P,\mathcal{C}})\cong P$.
Obviously, the labeling $\mathcal{C}$ is a weak coordinatization and it does not satisfy the the condition (C2).

\begin{theorem}\label{sect4-2t}
Let $P$ be a super-atomic lattice. Then the labeling $\mathcal{C}$ of $P$ as defined by (\ref{eq5.2}) is a coordinatization
if and only if for each $p\in (P\setminus\emph{atoms}(P))\setminus0$, either $N([a_i\vee a_k, 1])\leq N([a_r\vee a_k, 1])$ or $N([a_j\vee a_k, 1])\leq N([a_r\vee a_k, 1])$
for any $a_k, a_r \in \emph{supp}(p)$ where $\{a_i, a_j\} \in B_p$.
\end{theorem}
\proof
Let $\mathcal{C}$ be a coordinatization. Then there exists an isomorphism $g:P\rightarrow LCM(C_{P,\mathcal{C}})$ with $g(a)= x(a)$ for
each $a\in \mbox{atoms}(P)$. Suppose that $p\in (P\setminus\emph{atoms}(P))\setminus0$ and there exist $a_k, a_r \in \mbox{supp}(p)$ such that
$$N([a_i\vee a_k, 1])> N([a_r\vee a_k, 1])\mbox{ and }N([a_j\vee a_k, 1]) > N([a_r\vee a_k, 1])$$
where $\{a_i, a_j\} \in B_p$. Let $a_{k}^{n_{y}}$ be the highest power of $a_k$ dividing $x(a_{y})$ for any $a_y\in \mbox{atoms}(P)$.
Then similar to the proof of (\ref{eq5.1}), we have $n_r > n_i$ and $n_r > n_j$.
Thus $a_{k}^{n_{r}} \nmid \mbox{lcm}\{x(a_i), x(a_j)\}$, i.e., $x(a_r)\nmid \mbox{lcm}\{x(a_i), x(a_j)\}$.
Note that $g(p)=\mbox{lcm}\{x(a_i), x(a_j)\}$. Thus $x(a_r)\nmid g(p)$. However, $a_r\in \mbox{supp}(p)$ yields that $g(a_{r})= x(a_r)\mid g(p)$, a contradiction.

Conversely, suppose that for all $p\in (P\setminus\emph{atoms}(P))\setminus 0$, either $N([a_i\vee a_k, 1])\leq N([a_r\vee a_k, 1])$
or $N([a_j\vee a_k, 1])\leq N([a_r\vee a_k, 1])$
for any $a_k, a_r \in \mbox{supp}(p)$ where $\{a_i, a_j\} \in B_p$.

In what follows, we first prove that $\triangle(a)= x(a)$ for all $a\in \mbox{atoms}(P)$. The proof will be completed by two parts.

(E1). Let $p\in (P\setminus\emph{atoms}(P))\setminus0$ and $\{a_i, a_j\} \in B_p$.
Now, we prove that \begin{equation}\label{equ5.3}x(a_{s})\mid \mbox{lcm}\{x(a_{i}), x(a_{j})\} \mbox{ if } a_{s}\in\mbox{supp}(p)\setminus\{a_{i}, a_{j}\}.\end{equation}

As $P$ is super-atomic, \begin{equation}\label{11114}a_{i}\vee a_{s}< p\mbox{ and }a_{j}\vee a_{s}< p.\end{equation}
Let $a_t \in \mbox{atoms}(P)$ and $a_{t}^{n_{y}}$ be the highest power of $a_t$ dividing $x(a_{y})$ for any $a_y \in \mbox{atoms}(P)$. We claim \begin{equation}\label{11116}a_{t}^{n_{s}}\mid \mbox{lcm}\{x(a_{i}), x(a_{j})\}.\end{equation}

If $a_t=a_s$ then clearly $n_s=0$. This follows that (\ref{11116}) holds.

If $a_{t}\neq a_{s}$, then there are two cases.

Case (1*). Suppose $a_{t}\notin \mbox{supp}(p)$. Then $a_{i}\vee a_{j}\vee a_{t}\vee a_{s}=p\vee a_t> p$. Thus $a_{i}\vee a_{j}\neq a_{i}\vee a_{j}\vee a_{t}\vee a_{s}$,
$a_{i}\vee a_{s}\neq a_{i}\vee a_{j}\vee a_{t}\vee a_{s}$ and $a_{j}\vee a_{s}\neq a_{i}\vee a_{j}\vee a_{t}\vee a_{s}$ by (\ref{11114}). We claim
that $a_{s}\vee a_{t}\neq a_{i}\vee a_{j}\vee a_{t}\vee a_{s}$. Otherwise, $a_s\vee a_t=a_{i}\vee a_{j}\vee a_{t}$
since $a_{i}\vee a_{j}\vee a_{t}=a_{i}\vee a_{j}\vee a_{t}\vee a_{s}$, which together with $P$ is super-atomic yields $s=i$ or $s=j$,
a contradiction. Therefore, either $a_{i}\vee a_{j}\vee a_{t}\vee a_{s}=a_{i}\vee a_{t}$ or $a_{i}\vee a_{j}\vee a_{t}\vee a_{s}=a_{j}\vee a_{t}$.

Obviously, $a_{i}\vee a_{j}\vee a_{t}\vee a_{s}=a_{i}\vee a_{t}$ implies that $a_{s}\vee a_{t}< a_{i}\vee a_{j}\vee a_{t}\vee a_{s}=a_i\vee a_t$.
Thus $ N([a_{s}\vee a_{t}, 1])> N([a_{i}\vee a_{t}, 1])$. Similar to the proof of (\ref{eq5.1}), we have $n_s< n_i$,
it follows that $a_{t}^{n_{s}}\mid \mbox{lcm}\{x(a_{i}), x(a_{j})\}$.
Similarly, we can prove that $a_{t}^{n_{s}}\mid \mbox{lcm}\{x(a_{i}), x(a_{j})\}$ when $a_{i}\vee a_{j}\vee a_{t}\vee a_{s}=a_{j}\vee a_{t}$.
Therefore, \begin{equation*}a_{t}^{n_{s}}\mid \mbox{lcm}\{x(a_{i}), x(a_{j})\}\end{equation*} in the case of $a_{t}\notin \mbox{supp}(p)$.

Case (2*). Suppose $a_{t}\in \mbox{supp}(p)$. From the hypotheses, either $N([a_i\vee a_t, 1])\leq N([a_s\vee a_t, 1])$ or $N([a_j\vee a_t, 1])\leq N([a_s\vee a_t, 1])$.
In the first case, similar to the proof of (\ref{eq5.1}), we have $n_s\leq n_i$. Thus $a_{t}^{n_{s}}\mid \mbox{lcm}\{x(a_{i}), x(a_{j})\}$.
Similarly, we can prove $a_{t}^{n_{s}}\mid \mbox{lcm}\{x(a_{i}), x(a_{j})\}$ when $N([a_j\vee a_t, 1])\leq N([a_s\vee a_t, 1])$. Therefore, \begin{equation*}a_{t}^{n_{s}}\mid \mbox{lcm}\{x(a_{i}), x(a_{j})\}\end{equation*} in the case of $a_{t}\in \mbox{supp}(p)$.

Therefore, by Case (1*) and Case (2*), we know that (\ref{11116}) holds if $a_{t}\neq a_{s}$.

From the definition of $\mathcal{C}$, we have that if $x$ is a variable of $x(a_{s})$ then $x\in \mbox{atoms}(P)$.
Thus, by formula (\ref{11116})$$x(a_{s})\mid \mbox{lcm}\{x(a_{i}), x(a_{j})\} \mbox{ if } a_{s}\in\mbox{supp}(p)\setminus\{a_{i}, a_{j}\},$$
i.e., (\ref{equ5.3}) is true.

(E2). We shall prove that \begin{equation}\label{0000027}\triangle(a)=x(a)\end{equation} for each $a\in \mbox{atoms}(P)$.

Indeed, let $q\in P$ and $q\geq a$. We claim that \begin{equation}\label{0000026}x(a)\mid \mbox{lcm}\{x(r): r \in T\}\end{equation} for any $T\in B_q$.

If $q=a$ then clearly (\ref{0000026}) holds.

If $q> a$, then there exist $a_u, a_v \in \mbox{supp}(q)$ such that $a_u\vee a_v=q$. As $P$ is super-atomic, $a_u$, $a_v\in T$ for any $T\in B_q$. Using (\ref{equ5.3}),
we have that $x(c)\mid \mbox{lcm}\{x(a_u), x(a_v)\}$ for all $c\in \mbox{supp}(q)$. Note that $a\in \mbox{supp}(q)$. Thus $x(a)\mid \mbox{lcm}\{x(a_u), x(a_v)\}$.
Therefore, (\ref{0000026}) is true.

Formulas (\ref{0000026}) implies that $x(a)\mid \mbox{lcm}\{x(r): r \in T\}$ for any $T\in B_q$ if $q\geq a$. Thus $x(a)\mid\triangle(a)$ by (\ref{equ1.2}).
Note that $\triangle(a)\mid x(a)$. Therefore $\triangle(a)=x(a)$, i.e., (\ref{0000027}) holds.

$\quad$

In order to prove that $\mathcal{C}$ is a coordinatization. By Lemma \ref{sect3-he}, it suffices to prove that $\mathcal{C}$ is a weak coordinatization finally.

For $q\in P$, define $g: P\longrightarrow LCM(I_{P, \mathcal{C}})$ to be a map such that
\begin{equation}\label{11117}g(q)=\mbox{lcm}\{\triangle(w): w\in \mbox{supp}(q)\}.\end{equation}
Obviously, $g$ is meet-preserving, join-preserving and surjection by B, C and D in the proof of Theorem \ref{sect3-1t}. Thus we only need to prove that $g$ is injective.

Clearly, if $g(u)=g(v)$ and $0\in \{u, v\}$ then $u=0=v$.

Next, suppose that $g(u)=g(v)$ and $u, v\in P\setminus0$. We shall prove $u=v$.

Indeed, if $v\nleq u$ then $\mbox{supp}(v)\setminus\mbox{supp}(u)\neq \emptyset$.
Let $a_t \in \mbox{supp}(v)\setminus\mbox{supp}(u)$. There are two cases as below.

Case (k1). If $u\in \mbox{atoms}(P)$, then by statement ($\ast$) in the proof of Theorem \ref{sect4-1t}, $u \nmid \triangle(u)$ and $u \mid \triangle(a_t)$. Hence $u\mid g(v)$ and $u\nmid g(u)$, contrary to $g(u)=g(v)$.

Case (k2). If $u\in (P\setminus\mbox{atoms}(P))\setminus0$, then there exists $\{a_{i}, a_{j}\}\in B_u$. Thus \begin{equation}\label{0000028}g(u)=\mbox{lcm}\{\triangle(a_i), \triangle(a_j)\}.\end{equation}

Obviously, $u=a_i \vee a_j\neq a_t \vee a_i \vee a_j$ since $a_t \nleq u$.
Thus either $a_t \vee a_i = a_t \vee a_i \vee a_j$ or $a_t \vee a_j = a_t \vee a_i \vee a_j$. In the first case, we first note that $a_t \vee a_i> a_j\vee a_i$. Then $N([a_{j}\vee a_{i}, 1])> N([a_{t}\vee a_{i}, 1])$. Let $a_{i}^{n_{j}}$ be the highest power of $a_i$ dividing $x(a_{j})$
and $a_{i}^{n_{t}}$ be the highest power of $a_i$ dividing $x(a_{t})$. Similar to the proof of (\ref{eq5.1}), $n_{t}> n_{j}$. Thus $x(a_t)\nmid x(a_j)$.
Again, by statement ($\ast$) in the proof of Theorem \ref{sect4-1t}, $a_i \nmid x(a_i)$ since $x(a_i)= \triangle(a_i)$, and this means that $x(a_t)\nmid x(a_i)$. Therefore, $x(a_t)\nmid \mbox{lcm}\{x(a_i), x(a_j)\}$. As $\triangle(a)=x(a)$ for any $a\in \mbox{atoms}(P)$, we have $\triangle(a_t)\nmid \mbox{lcm}\{\triangle(a_i), \triangle(a_j)\}$. From formulas (\ref{11117}) and (\ref{0000028}),
we have $\triangle(a_t)\nmid g(u)$ but $\triangle(a_t)\mid g(v)$ since $a_t \in \mbox{supp}(v)$, contrary to $g(u)=g(v)$. In the second case, with analogous proof to the first case of $a_t \vee a_i = a_t \vee a_i \vee a_j$, one can deduce a contradiction.

Cases (k1) and (k2) tell us that the assumption of $v\nleq u$ will yield a contradiction. Hence $v\leq u$.

Arguing as above, we can prove that $u\leq v$. Therefore, $u=v$.

Consequently, $g$ is injective.
\endproof

Using Theorem \ref{sect4-2t}, we can determine whether the labeling, defined by (\ref{eq5.2}), of a super-atomic lattice is a coordinatization.

As a conclusion of this section, we shall consider when the labeling, defined by (\ref{eq5.2}), of a non-super-atomic lattice is also a coordinatization.

\begin{lemma}\label{sect5-2t} Let $P, Q\in\mathcal{L}(n)$ with $\emph{atoms}(P)=\emph{atoms}(Q)=\{1, 2, \cdots, n\}$. If $\mathcal{S}_P\setminus\mathcal{S}_Q=\{S\}$ then $S$ is meet-irreducible in $(\mathcal{S}_P, \subseteq)$.
\end{lemma}
\proof If $S$ is not a meet-irreducible in $(\mathcal{S}_P, \subseteq)$, then there exist two different elements $S_1, S_2\in \mathcal{S}_P$ such that $S_1\succ S$
and $S_2\succ S$ in lattice $(\mathcal{S}_P, \subseteq)$. Note that $S_1, S_2\in \mathcal{S}_Q$. We claim that $\bigvee_{t\in S}\{t\} = S_1$ in
lattice $(\mathcal{S}_Q, \subseteq)$.
Otherwise, we have $\bigvee_{t\in S}\{t\} = R\subsetneq S_1$ for some $R\in \mathcal{S}_Q$ in lattice $(\mathcal{S}_Q, \subseteq)$. Clearly, $S\subseteq R\subsetneq S_1$.
As $S\notin \mathcal{S}_Q$, $S\neq R$ which means that $S\subsetneq R$. Therefore, $S\subsetneq R\subsetneq S_1$, which together with $S, R, S_1 \in \mathcal{S}_P$
yields that $S_1 \nsucc S$ in lattice $(\mathcal{S}_P, \subseteq)$, a contradiction. Consequently, $\bigvee_{t\in S}\{t\} = S_1$ in $(\mathcal{S}_Q, \subseteq)$.
Similarly, we also have $\bigvee_{t\in S}\{t\} = S_2$ in $(\mathcal{S}_Q, \subseteq)$. Therefore, $S_1=S_2$, contrary to $S_1\neq S_2$.
\endproof

Let $P\in \mathcal{L}(n)$ with $\mbox{atoms}(P)=\{a_1, a_2,\cdots, a_n\}$. Next we denote by $\mathcal{C}_P$ the labeling of $P$
defined by (\ref{eq5.2}), that is, $m_c=\prod_{a_i\in \mbox{supp}(c)}a_i$ for any $c\in P\setminus0$.
Note that $(\mathcal{S}_P, \subseteq)$ is the lattice corresponding to $P$ (see Section 2).
Then for any $C\in \mathcal{S}_P\setminus\emptyset$, we have that $m_C=\prod_{a_i\in C}a_i$ where $C$ corresponds to $c$.
Again, we denote by $x_{P}(\{a_i\})$ the monomials corresponding to $(\mathcal{S}_P, \subseteq)$ defined by (\ref{eqhe0}).
Then we define $C_{\mathcal{S}_P, \mathcal{C}_P}$ as the ideal generated by monomials $x_{P}(\{a_i\})$ for each $i\in\{1, 2, \cdots, n\}$.
We denote by $\triangle_{P}(\{a_i\})$ the monomials corresponding to $(\mathcal{S}_P, \subseteq)$ defined by (\ref{equ1.2}),
and define $I_{\mathcal{S}_P, \mathcal{C}_P}$ as the ideal generated by monomials $\triangle_{P}(\{a_i\})$ for each $i\in\{1, 2, \cdots, n\}$.
Then we have the following theorem.

\begin{theorem}\label{sect5-3t}
Let $(\mathcal{S}_Q, \subseteq), (\mathcal{S}_P, \subseteq), (\mathcal{S}_R, \subseteq) \in \mathcal{L}(n)$ and $(\mathcal{S}_R, \subseteq)$ be a super-atomic lattice.
If $\mathcal{S}_P\subseteq \mathcal{S}_R$, $\mathcal{S}_P\setminus \mathcal{S}_Q=\{S\}$ and $\mathcal{C}_P$ be a coordinatization,
then $\mathcal{C}_{Q}$ is a coordinatization if and only if $\triangle_{Q}(\{a_k\})= x_{Q}(\{a_k\})$
for any $k\in\{1, 2, \cdots, n\}$.
\end{theorem}
\proof
We only need to show the sufficiency of theorem since the necessity is obvious.
We first note that $I_{\mathcal{S}_Q, \mathcal{C}_Q} = C_{\mathcal{S}_Q, \mathcal{C}_Q}$ since $\triangle_{Q}(\{a_k\})= x_{Q}(\{a_k\})$
for any $k\in\{1, 2, \cdots, n\}$.
Define a map $h: (\mathcal{S}_Q, \subseteq) \rightarrow LCM(I_{\mathcal{S}_Q, \mathcal{C}_Q})=LCM(C_{\mathcal{S}_Q, \mathcal{C}_Q})$
as $$h(C)= \mbox{lcm}\{\triangle_{Q}(\{a_i\}): a_i\in C\}= \mbox{lcm}\{x_{Q}(\{a_i\}): a_i\in C\}$$
for any $C\in \mathcal{S}_Q$. According to Lemma \ref{sect3-he}, we just need to prove that $\mathcal{C}_{Q}$ is a weak coordinatization, i.e.,
we just need to prove $h$ is an isomorphism. By B, C and D in the proof of Theorem \ref{sect3-1t}, one can check that $h$ is meet-preserving, join-preserving and surjective.
Now, we shall prove that $h$ is injective.

For $C\in \mathcal{S}_P$, we define a map $g: (\mathcal{S}_P,\subseteq) \rightarrow LCM(C_{\mathcal{S}_P, \mathcal{C}_P})$ such that $$g(C)= \mbox{lcm}\{x_{P}(\{a_i\}): a_i\in C\}.$$
Obviously, $g$ is an isomorphism from $(\mathcal{S}_P, \subseteq)$ to $LCM(C_{\mathcal{S}_P, \mathcal{C}_P})$ since $\mathcal{C}_P$ is a coordinatization.

By Lemma \ref{sect5-2t} there exists exactly one element $T\in \mathcal{S}_P$ such that $T\succ S$ in lattice $(\mathcal{S}_P, \subseteq)$.
Clearly $S\notin \mbox{atoms}(\mathcal{S}_P)\bigcup \{\emptyset\}$.

If $a_j\in \{a_1,a_2,\cdots,a_n\}\setminus S$ then $S \notin \lceil\{a_j\}\rceil_{P}$ since $\{a_j\}\nsubseteq S$.
Thus $\lceil\{a_j\}\rceil_{P}=\lceil\{a_j\}\rceil_{Q}$, which implies that
\begin{equation}\label{I1}x_{Q}(\{a_j\})=\prod_{C\in \lceil\{a_j\}\rceil_{Q}^{c}}m_C=\prod_{C\in \mathcal{S}_Q\setminus\lceil\{a_j\}\rceil_{Q}}m_C=\frac{\prod_{C\in \mathcal{S}_P\setminus\lceil\{a_j\}\rceil_{P}}m_C}{m_S}=
\frac{x_{P}(\{a_j\})}{\prod_{a_i\in S} a_{i}}.\end{equation}

If $a_j\in S$ then $S\in \lceil\{a_j\}\rceil_{P}$ since $\{a_j\}\subseteq S$. Thus $\lceil\{a_j\}\rceil_{P}=\lceil\{a_j\}\rceil_{Q}\cup \{S\}$, which implies that
\begin{equation}\label{I2}x_{Q}(\{a_j\})=\prod_{C\in \lceil\{a_j\}\rceil_{Q}^{c}}m_C=\prod_{C\in \mathcal{S}_Q\setminus\lceil\{a_j\}\rceil_{Q}}m_C= \prod_{C\in \mathcal{S}_P\setminus\lceil\{a_j\}\rceil_{P}}m_C = x_{P}(\{a_j\}).\end{equation}

The following proof is completed by three parts.

(I) Let $C_1, D_1 \in \mathcal{S}_Q$. If $h(C_1)=h(D_1)$ and $C_{1}\subseteq D_1$ then $C_1= D_1$.

Suppose that $C_1\neq D_1$. Then $C_1\subsetneq D_1$. Thus there exists $C_2\in \mathcal{S}_Q$ such that
\begin{equation}\label{chen4}C_1\prec C_2\subseteq D_1\mbox{ in }(\mathcal{S}_Q, \subseteq),\end{equation}
and \begin{equation}\label{0000029}h(C_1)=h(C_2)\end{equation} since $h$ is meet-preserving.

Clearly, if $C_1=\emptyset$ then $h(C_1)=1=h(D_1)$, and which implies that $C_1=D_1$. Next, we suppose that $C_1\in \mathcal{S}_Q\setminus\emptyset$.

If $C_1 \in \mbox{atoms}((\mathcal{S}_Q, \subseteq))$, then let $C_1=\{a_u\}$. Clearly, there exists $\{a_v\}\subseteq C_2$ such that $\{a_v\}\neq \{a_u\}$ by (\ref{chen4}). By statement ($\ast$), we know that $a_u \nmid \triangle_{Q}(\{a_u\})$ and $a_u \mid \triangle_{Q}(\{a_v\})$.  Hence, $a_u\mid h(C_2)$ and $a_u\nmid h(C_1)$, contrary to formula (\ref{0000029}).

If $C_1 \in (\mathcal{S}_Q\setminus \mbox{atoms}((\mathcal{S}_Q, \subseteq)))\setminus\emptyset$, then there
exist $\{a_i\}, \{a_j\}\in \mbox{atoms}((\mathcal{S}_Q, \subseteq))$
such that \begin{equation}\label{0000030}C_1=\{a_i\}\vee \{a_j\}\end{equation} in $(\mathcal{S}_Q, \subseteq)$
since $(\mathcal{S}_R, \subseteq)$ is super-atomic and $\mathcal{S}_Q \subseteq \mathcal{S}_R$.
Further, by (\ref{chen4}), there exists $\{a_k\}\in \mbox{atoms}((\mathcal{S}_Q, \subseteq))$
such that \begin{equation}\label{0000031}C_2=\{a_i\}\vee \{a_j\}\vee \{a_k\}\end{equation}
in $(\mathcal{S}_Q, \subseteq)$. Using formulas (\ref{0000029}), (\ref{0000030}) and (\ref{0000031}),
we have
\begin{equation}\label{0000032}h(C_1)=\mbox{lcm}\{x_{Q}(\{a_i\}), x_{Q}(\{a_j\})\}=\mbox{lcm}\{x_{Q}(\{a_i\}), x_{Q}(\{a_j\}), x_{Q}(\{a_k\})\}=h(C_2).\end{equation}
Thus we shall distinguish the six types as follows.

In what follows, for convenience, let $a_{y}^{m_{x_y}}$ be the highest power of $a_{y}$ dividing $x_{P}(\{a_x\})$
and $a_{y}^{n_{x_y}}$ be the highest power of $a_{y}$ dividing $x_{Q}(\{a_x\})$ for any $x, y\in \{1, 2, \cdots, n\}$.

$\quad$\\
Type 1. $a_i, a_j, a_k \in S$.

We first claim that $C_1\neq T$. If $C_1=T$, then $\{a_i\}\vee \{a_j\}=S$ in $(\mathcal{S}_P, \subseteq)$ since $a_i, a_j \in S$.
Thus $\{a_i\}\vee \{a_j\}\vee \{a_k\}= S$ in $(\mathcal{S}_P, \subseteq)$ since $a_k\in S$. By formula (\ref{0000031}), $C_2=T$.
So that $C_2=C_1$, a contradiction. Hence $C_1\neq T$, and $C_1\subsetneq T$ since $\{a_i,a_j\}\subseteq T$.
Therefore,  \begin{equation}\label{0000033}\{a_i\}\vee \{a_j\}= C_1\subsetneq S\end{equation} in $(\mathcal{S}_P, \subseteq)$.

Using formula (\ref{I2}), we have $x_{Q}(\{a_t\})=x_{P}(\{a_t\})$ for any $t\in \{i, j, k\}$. Then \begin{equation}\label{0000034}\mbox{lcm}\{x_{P}(\{a_i\}), x_{P}(\{a_j\})\}=\mbox{lcm}\{x_{P}(\{a_i\}), x_{P}(\{a_j\}), x_{P}(\{a_k\})\}\end{equation} by formula (\ref{0000032}). There are two subcases as below.

Subcase (li). If $C_2=T$ then $\{a_i\}\vee \{a_j\}\vee \{a_k\}= S$ in $(\mathcal{S}_P, \subseteq)$ since $a_i,a_j,a_k\in S$, which together with formulas (\ref{0000033}) and (\ref{0000034}) implies that $g(C_1)=g(S)$. However, $g(C_1) < g(S)$ since $C_1\subsetneq S$ and $g$ is isomorphic, a contradiction.

Subcase (2i). If $C_2\neq T$ then $\{a_i\}\vee \{a_j\}\vee \{a_k\}= C_2$ in $(\mathcal{S}_P, \subseteq)$. By formulas (\ref{0000033})
and (\ref{0000034}), $g(C_1)=g(C_2)$, contrary to $g(C_1)< g(C_2)$.

$\quad$\\
Type 2. $a_i, a_j, a_k\notin S$.

By formula (\ref{I1}), $x_{P}(\{a_t\})=(\prod_{a_r\in S}a_{r})\ast x_{Q}(\{a_t\})$ for any $t\in \{i,j,k\}$.
Then $h(C_1)= h(C_2)$ implies that $h(C_1)\ast \prod_{a_r\in S}a_{r}= h(C_2)\ast \prod_{a_r\in S}a_{r}$.
Further, by formula (\ref{0000032}),
$$\mbox{lcm}\{x_{P}(\{a_i\}), x_{P}(\{a_j\})\}=\mbox{lcm}\{x_{P}(\{a_i\}), x_{P}(\{a_j\}), x_{P}(\{a_k\})\}.$$
On the other hand, as $a_i, a_j, a_k\notin S$, $\{a_i\}\vee \{a_j\}= C_1$ and $\{a_i\}\vee \{a_j\}\vee \{a_k\}= C_2$ in $(\mathcal{S}_P, \subseteq)$
obviously. Therefore, $g(C_1)= g(C_2)$, contrary to $g(C_1)< g(C_2)$.

$\quad$\\
Type 3. $a_i, a_j\notin S$ and $a_k \in S$.

By formulas (\ref{I2}) and (\ref{0000032}), we have that
$$\mbox{lcm}\{x_{Q}(\{a_i\}), x_{Q}(\{a_j\})\}=\mbox{lcm}\{x_{Q}(\{a_i\}), x_{Q}(\{a_j\}), x_{P}(\{a_k\})\}.$$
Thus $x_{P}(\{a_k\})\mid \mbox{lcm}\{x_{Q}(\{a_i\}), x_{Q}(\{a_j\})\}$. Similar to the proof of Type 2, we know that  $\{a_i\}\vee \{a_j\}= C_1$, $\{a_i\}\vee \{a_j\}\vee \{a_k\}= C_2$ in $(\mathcal{S}_P, \subseteq)$ and $x_{P}(\{a_t\})=(\prod_{a_r\in S}a_{r})\ast x_{Q}(\{a_t\})$ for any $t\in \{i,j\}$. Thus $x_{P}(\{a_k\})\mid \mbox{lcm}\{x_{P}(\{a_i\}), x_{P}(\{a_j\})\}$, which implies that $$\mbox{lcm}\{x_{P}(\{a_i\}), x_{P}(\{a_j\})\}=\mbox{lcm}\{x_{P}(\{a_i\}), x_{P}(\{a_j\}), x_{P}(\{a_k\})\}.$$
Therefore, $g(C_1)=g(C_2)$, contrary to $g(C_1)< g(C_2)$.

$\quad$\\
Type 4. $a_i\in S$, $a_j\notin S$ and $a_k\in S$.

Using (\ref{I2}) and (\ref{0000032}), we have that
$$\mbox{lcm}\{x_{P}(\{a_i\}), x_{Q}(\{a_j\})\}=\mbox{lcm}\{x_{P}(\{a_i\}), x_{Q}(\{a_j\}), x_{P}(\{a_k\})\}.$$
Similar to the proof of Type 3, we have that $x_{P}(\{a_k\})\mid \mbox{lcm}\{x_{P}(\{a_i\}), x_{P}(\{a_j\})\}$ and $g(C_1)=g(C_2)$ with $\{a_i\}\vee \{a_j\}= C_1$ and $\{a_i\}\vee \{a_j\}\vee \{a_k\}= C_2$ in $(\mathcal{S}_P, \subseteq)$, contrary to $g(C_1)< g(C_2)$.

$\quad$\\
Type 5. $a_i,a_j\in S$ and $a_k\notin S$.

Using (\ref{I2}) and (\ref{0000032}), we have that
$$\mbox{lcm}\{x_{P}(\{a_i\}), x_{P}(\{a_j\})\}=\mbox{lcm}\{x_{P}(\{a_i\}), x_{P}(\{a_j\}), x_{Q}(\{a_k\})\}.$$
Then \begin{equation}\label{0eq1}x_{Q}(\{a_k\})\mid \mbox{lcm}\{x_{P}(\{a_i\}), x_{P}(\{a_j\})\}.\end{equation}

Using (\ref{I1}), we have $x_{P}(\{a_k\})=(\prod_{a_r\in S}a_{r})\ast x_{Q}(\{a_k\})$.
Thus $n_{k_{i}}+1 =m_{k_{i}}$ since $a_i\in S$.  We note that $\{a_i\}\vee \{a_j\}\vee \{a_k\}= C_2$ in $(\mathcal{S}_P, \subseteq)$
since $a_k\notin S$. Then $$\{a_i\}\vee \{a_k\}=C_2\mbox{ or }\{a_j\}\vee \{a_k\}=C_2$$ in $(\mathcal{S}_P, \subseteq)$ since
$\mathcal{S}_P \subseteq \mathcal{S}_R$ and $(\mathcal{S}_R, \subseteq)$ is super-atomic. There are two subcases.

Subcase 1. If $C_1=T$ then $\{a_i\}\vee \{a_j\}= S$ in $(\mathcal{S}_P, \subseteq)$. Thus $S \subsetneq C_1=T\subsetneq C_2$.

Assume that $\{a_i\}\vee \{a_k\}= C_2$ in $(\mathcal{S}_P, \subseteq)$.
Then we have that $$N([\{a_i\}\vee \{a_j\}, 1])\geq N([\{a_i\}\vee \{a_k\}, 1])+2$$
in $(\mathcal{S}_P, \subseteq)$ since $S \subsetneq T\subsetneq C_2$.
Similar to the proof of (\ref{eq5.1}), we have that $m_{k_{i}} \geq m_{j_{i}}+2$. Thus $n_{k_{i}} \geq m_{j_{i}}+1$ which implies that $x_{Q}(\{a_k\})\nmid x_{P}(\{a_j\})$.
From Lemma \ref{sect3-he}, $\triangle_{P}(\{a_i\})= x_{P}(\{a_i\})$ since $\mathcal{C}_P$ is a coordinatization. Further, by statement ($\ast$), we know that $a_i\nmid x_{P}(\{a_i\})$. Therefore, $x_{Q}(\{a_k\})\nmid \mbox{lcm}\{x_{P}(\{a_i\}), x_{P}(\{a_j\})\}$, contrary to (\ref{0eq1}).

If $\{a_j\}\vee \{a_k\}= C_2$ in $(\mathcal{S}_P, \subseteq)$, then with analogous proof to
the case of $\{a_i\}\vee \{a_k\}= C_2$ in $(\mathcal{S}_P, \subseteq)$ one may get a contradiction.

Subcase 2. If $C_1\neq T$ then $C_1 \subsetneq S$ and $\{a_i\}\vee \{a_j\}= C_1$ in $(\mathcal{S}_P, \subseteq)$ by the proof of Type 1.

Suppose that $\{a_i\}\vee \{a_k\}= C_2$ in $(\mathcal{S}_P, \subseteq)$. Then $$N([\{a_i\}\vee \{a_j\}, 1])> N([\{a_i\}\vee \{a_k\}, 1])$$
in $(\mathcal{S}_P, \subseteq)$ since $C_1\subsetneq C_2$. Note that $C_2\nsubseteq S$ since $a_k\notin S$.
Thus $$N([\{a_i\}\vee \{a_j\}, 1])\geq N([\{a_i\}\vee \{a_k\}, 1]) +2$$ in $(\mathcal{S}_P, \subseteq)$ since $C_1\subsetneq S$.

Similar to Subcase 1, one can prove that $x_{Q}(\{a_k\})\nmid \mbox{lcm}\{x_{P}(\{a_i\}), x_{P}(\{a_j\})\}$, contrary to (\ref{0eq1}).

If $\{a_j\}\vee \{a_k\}= C_2$ in $(\mathcal{S}_P, \subseteq)$, then with analogous proof to
the case of $\{a_i\}\vee \{a_k\}= C_2$ in $(\mathcal{S}_P, \subseteq)$ one may get a contradiction.

$\quad$\\
Type 6. $a_i\in S$ and $a_j,a_ k \notin S$.

By (\ref{I2}) and (\ref{0000032}), we have that
$$\mbox{lcm}\{x_{P}(\{a_i\}), x_{Q}(\{a_j\})\}=\mbox{lcm}\{x_{P}(\{a_i\}), x_{Q}(\{a_j\}), x_{Q}(\{a_k\})\}.$$
Thus \begin{equation}\label{1eq1}x_{Q}(\{a_k\})\mid \mbox{lcm}\{x_{P}(\{a_i\}), x_{Q}(\{a_j\})\}.\end{equation}
Clearly, $\{a_i\}\vee \{a_j\}=C_1$ and $\{a_i\}\vee \{a_j\}\vee \{a_k\}=C_2$ in $(\mathcal{S}_P, \subseteq)$
since $a_j, a_k \notin S$. By the proof of Type 5, we know that
$$\{a_i\}\vee \{a_k\}=C_2\mbox{ or }\{a_j\}\vee \{a_k\}=C_2$$ in $(\mathcal{S}_P, \subseteq)$. There are two subcases.

Subcase (i). If $\{a_i\}\vee \{a_k\}=C_2$ in $(\mathcal{S}_P, \subseteq)$. By the proof of Type 5, we have that $$N([\{a_i\}\vee \{a_j\}, 1])> N([\{a_i\}\vee \{a_k\}, 1])$$ in $(\mathcal{S}_P, \subseteq)$. Clearly, $m_{k_{i}}> m_{j_{i}}$, i.e., $x_{P}(\{a_k\})\nmid x_{P}(\{a_j\})$. Using (\ref{I1}), we have $$x_{P}(\{a_j\})=(\prod_{a_r\in S}a_{r})\ast x_{Q}(\{a_j\})\mbox{ and }x_{P}(\{a_k\})=(\prod_{a_r\in S}a_{r})\ast x_{Q}(\{a_k\}).$$ Hence $x_{Q}(\{a_k\})\nmid x_{Q}(\{a_j\})$.

From Lemma \ref{sect3-he}, $\triangle_{P}(\{a_i\})= x_{P}(\{a_i\})$ since $\mathcal{C}_P$ is a coordinatization. Further, by statement ($\ast$), $a_i\nmid x_{P}(\{a_i\})$. Thus $x_{Q}(\{a_k\})\nmid \mbox{lcm}\{x_{P}(\{a_i\}), x_{Q}(\{a_j\})\}$, contrary to the formula (\ref{1eq1}).

Subcase (ii). If $\{a_j\}\vee \{a_k\}=C_2$ in $(\mathcal{S}_P, \subseteq)$.
We note that $$N([\{a_i\}\vee \{a_j\}, 1])> N([\{a_j\}\vee \{a_k\}, 1])$$ in $(\mathcal{S}_Q, \subseteq)$.
Clearly, $n_{k_{j}}> n_{i_{j}}$. Again, we know that $n_{i_{j}}=m_{i_{j}}$ since $x_{P}(\{a_i\})=x_{Q}(\{a_i\})$. Hence $x_{Q}(\{a_k\})\nmid x_{P}(\{a_i\})$.  Since $\triangle_{Q}(\{a_j\})= x_{Q}(\{a_j\})$, we have $a_j\nmid x_{Q}(\{a_j\})$ by statement ($\ast$). Therefore, $x_{Q}(\{a_k\})\nmid \mbox{lcm}\{x_{P}(\{a_i\}), x_{Q}(\{a_j\})\}$ , contrary to the formula (\ref{1eq1}).

Types 1-6 tell us that if $h(C_1)=h(D_1)$ and $C_{1}\subseteq D_1$ then $C_1= D_1$.

$\quad$

Similar to (I), we can prove that

(II) If $h(C_1)= h(D_1)$ and $C_1\supseteq D_1$ then $C_1=D_1$.

$\quad$

(III) If $h(C_1)= h(D_1)$ then $C_1\subseteq D_1$ or $C_1\supseteq D_1$.

Assume that $C_1\| D_1$. Let $\{a_i\}\vee \{a_j\}=C_1$ and $\{a_k\}\vee \{a_e\}=D_1$ in $(\mathcal{S}_Q, \subseteq)$.
Then $ C= C_1\vee D_1 = \{a_i\}\vee \{a_j\}\vee \{a_k\} \vee \{a_e\}\supsetneq\{a_i\}\vee \{a_j\}= C_1$ in $(\mathcal{S}_Q, \subseteq)$.
Thus by (I), we have that $h(C_1)< h(C)$. This follows that $$\mbox{lcm}\{x_{Q}(\{a_i\}), x_{Q}(\{a_j\})\}< \mbox{lcm}\{x_{Q}(\{a_i\}), x_{Q}(\{a_j\}), x_{Q}(\{a_k\}), x_{Q}(\{a_e\})\}.$$
Therefore,
\begin{equation}\label{0eq2}x_{Q}(\{a_k\})\nmid \mbox{lcm}\{x_{Q}(\{a_i\}), x_{Q}(\{a_j\})\} \mbox{ or } x_{Q}(\{a_e\})\nmid \mbox{lcm}\{x_{Q}(\{a_i\}), x_{Q}(\{a_j\})\},\end{equation} and formula (\ref{0eq2}) imply that $$h(D_1)=\mbox{lcm}\{x_{Q}(\{a_k\}), x_{Q}(\{a_e\})\} \nmid  \mbox{lcm}\{x_{Q}(\{a_i\}), x_{Q}(\{a_j\})\}= h(C_1),$$ i.e., $h(C_1)\neq h(D_1)$, a contradiction.

From (I), (II) and (III), we know that the map $h$ is injective.
\endproof

The following example will illustrate Theorem \ref{sect5-3t}.
\begin{example}\label{sect5-2e}
\emph{Let
$$\mathcal{S}_P=\{\{a_1,a_2,a_3,a_4\}, \{a_2,a_3, a_4\}, \{a_1,a_3,a_4\},\{a_3,a_4\},\{a_2,a_3\},\{a_1,a_4\},\{a_1\}, \{a_2\}, \{a_3\}, $$
$\{a_4\},\emptyset\}.$
It is easy to see that $(\mathcal{S}_P,\subseteq)$ is a supper-atomic lattice in $\mathcal{L}(4)$. Denote $\mathcal{C}_P$ as a labeling of $\mathcal{S}_P$ defined by (\ref{eq5.2}). Then $C_{\mathcal{S}_P,\mathcal{C}_{P}}=\{a_{2}^{3}a_{3}^{4}a_{4}^{3}, a_{1}^{3}a_{3}^{3}a_{4}^{4}, a_{1}^{2}a_{2}a_{4}^{2}, a_{1}a_{2}^{2}a_{3}^{2}\}$. Clearly, the labeling $\mathcal{C}_{P}$ is a coordinatization.}

\emph{Let $\mathcal{S}_Q=\mathcal{S}_P\setminus\{\{a_2,a_3,a_4\}\}$.
Clearly $x_{Q}(\{a_i\})=\triangle_{Q}(\{a_i\})$ for any $i\in \{1,2,3,4\}$. Then by Theorem \ref{sect5-3t},
$$C_{\mathcal{S}_Q,\mathcal{C}_{Q}}=\{a_{2}^{2}a_{3}^{3}a_{4}^{2}, a_{1}^{3}a_{3}^{3}a_{4}^{4}, a_{1}^{2}a_{2}a_{4}^{2}, a_{1}a_{2}^{2}a_{3}^{2}\}.$$ Further, one can check that $LCM(C_{\mathcal{S}_Q,\mathcal{C}_{Q}})\cong (\mathcal{S}_Q, \subseteq)$, i.e., $\mathcal{C}_{Q}$ is a coordinatization.}
\end{example}

\section{Conclusions}
This paper studied monomial ideals by their associated lcm-lattices. It first introduced notions of weak coordinatizations which have weaker hypotheses than coordinatizations, and showed the characterizations of all such weak coordinatizations which partly answer the problem arisen by Mapes in \cite{Mapes13}. It then defined a finite super-atomic lattice in $\mathcal{L}(n)$ which are used to investigate the structures of $\mathcal{L}(n)$ and to identify a specific labeling, given by us, of finite atomic lattice is the weak coordinatizations. It will be very interesting to study a minimal free resolution of $R/M$ by our results in the future.
\section*{Acknowledgments}
The authors thank the referees for their valuable comments and suggestions.

\end{document}